\documentclass[12pt]{article}
\usepackage{amssymb}
\usepackage{epsf}
\usepackage{amsmath, amsthm}
\usepackage{version}
\usepackage{color}
\usepackage{graphics}
\numberwithin{equation}{section}
\newtheorem{theorem}[equation]{Theorem}
\newtheorem{corollary}[equation]{Corollary}
\newtheorem{definition}[equation]{Definition}
\newtheorem{proposition}[equation]{Proposition}
\newtheorem{lemma}[equation]{Lemma}
\newtheorem{remark}[equation]{Remark}

\makeatletter
\newcommand{\customlabel}[2]{%
   \protected@edef\@currentlabel{#1}%
   \label{#2}%
}
\makeatother
\newtheorem*{newremark}{Remark}

\newcommand\beq{\begin{equation}}
\newcommand\ds{\displaystyle}
\newcommand\eeq{\end{equation}}
\newcommand\re{\mathrm {Re~}}
\renewcommand\Re{\mathrm {Re~}}
\newcommand\im{\mathrm {Im~}}

\newcommand\nt{\stackrel{\rm  nt }{\to}}
\newcommand\al{\alpha}
\newcommand\de{\delta}
\newcommand\ga{\gamma}
\newcommand\ep{\varepsilon}
\newcommand\la{\lambda}
\renewcommand\ln{{\lambda_n}}

\newcommand{\ph}{\varphi}
\newcommand\si{\sigma}

\newcommand\D{\mathbb D}
\newcommand\T{\mathbb T}
\newcommand\C{\mathbb C}
\newcommand\M{\mathcal{M}}
\newcommand\R{\mathbb R}
\newcommand\Ha{\mathbb{H}}
\newcommand{\HH}{\Ha\times\Ha}
\newcommand\Htau{\mathbb{H}(\tau)}

\newcommand\Pick{\mathcal P}
\newcommand\Schur{\mathcal{S}}

\newcommand\lan{\langle}
\newcommand\ran{\rangle}

\newcommand\half{{\tfrac{1}{2}}}

\newcommand\df{\stackrel{\rm def}{=}}
\newcommand\blue{\color{blue}}
\newcommand\black{\color{black}}

\newcommand\nn{\nonumber}
\newcommand\tb{\partial (\D^2)} 
\renewcommand\phi{\varphi}
\renewcommand\epsilon{\varepsilon}

\def\cjw{Julia-Carath\'eodory\ }

\numberwithin{equation}{section}
\title{A Carath\'eodory theorem for the bidisk via Hilbert space methods}
\author{Jim Agler
\thanks{Partially supported by National Science Foundation Grant
DMS 0801259}
\\ U.C. San Diego\\
\and
John E. M\raise.5ex\hbox{c}Carthy
\thanks{Partially supported by National Science Foundation Grant DMS 0501079}
\\ Washington University, St. Louis \\
\and
N. J. Young
\thanks{Supported by EPSRC Grant EP/G000018/1}
\\  Newcastle University}

\begin{document}
\bibliographystyle{plain}

\maketitle
\begin{abstract}
If $\ph$ is an analytic function bounded by $1$ on the bidisk $\D^2$ and $\tau\in\tb$ is a point at which $\ph$ has an angular gradient $\nabla\ph(\tau)$ then $\nabla\ph(\la) \to \nabla\ph(\tau)$ as $\la\to\tau$ nontangentially in $\D^2$.  This is an analog for the bidisk of a classical theorem of Carath\'eodory for the disk.

For $\ph$ as above, if $\tau\in\tb$ is such that the $\liminf$ of $(1-|\ph(\la)|)/(1-\|\la\|)$ as $\la\to\tau$ is finite then the directional derivative $D_{-\de}\ph(\tau)$ exists for all appropriate directions $\de\in\C^2$.  Moreover, one can associate with $\ph$ and $\tau$ an analytic function $h$ in the Pick class such that the value of the directional derivative can be expressed in terms of $h$.
\end{abstract}

\section{Introduction} \label{intro}

\begin{figure}
\begin{center}
\begin{picture}(200,200)(,)

\epsfxsize=200pt
\epsffile{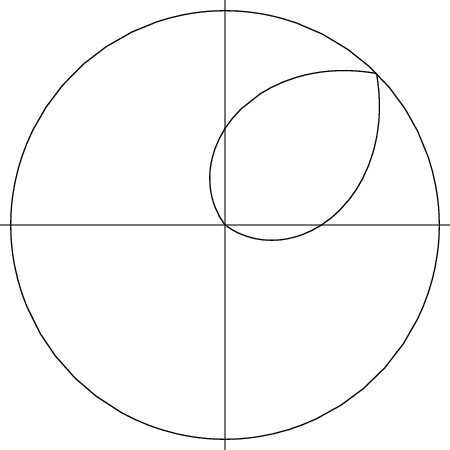}

\end{picture}
\caption{A non-tangential approach region at $\tau = e^{i\pi / 4}$}
\label{piccb1}
\end{center}
\end{figure}

To what extent do classical theorems on analytic functions in the unit disk $\D$ have analogs for the polydisk?    We are interested in results on the
 boundary behavior of  analytic functions
that map  $\D$ to $\D^-$; this set of functions is called
 the Schur class. (We shall use $S^-$ to denote the closure of any set $S$).
 The Julia-Carath\'eodory Theorem, due to G. Julia \cite{ju20} in 1920 and C. Carath\'eodory \cite{car29} in 1929, asserts that if a function $\ph$ in the Schur class satisfies the weak regularity condition,
\beq\label{A}
\liminf_{\la\to\tau} \frac{1-|\ph(\la)|}{1-|\la|} < \infty
\eeq
 at a point $\tau$ on the unit circle $\T$, then it also satisfies three further, apparently stronger, conditions (B to D of Theorem \ref{thma1} below). These state:

(B) The quotient on the right-hand side of (\ref{A}) has a non-tangential limit as $\la$ tends to $\tau$. (Non-tangential, which we shall define precisely in Section \ref{nontang},
means that $\la$ approaches $\tau$ from within a region like the one in Figure \ref{piccb1}.)

(C) The function $\ph$ has both a non-tangential limit $\omega \in \T$ at $\tau$ and also an angular derivative $\eta \in \C$, that is the difference quotient
\[
\frac{\ph(\la) - \omega}{\la - \tau}
\]
has a non-tangential limit $\eta$ at $\tau$.

(D) There exist $\omega$ in $\T$ and $\eta$ in $\C$ so that at $\tau$, $\ph(\la)$ tends to $\omega$ non-tangentially 
and $\phi'(\la)$ tends to $\eta$ non-tangentially.

 Futhermore, if (\ref{A}) holds, then a form (\ref{jul7}) of boundary Schwarz-Pick inequality holds.

It is known that some parts of the Julia-Carath\'eodory Theorem do extend to the Schur class of the  polydisk, but that some do not. The Schur
class of the polydisk $\D^n$ is the set of analytic functions from $\D^n$ to $\D^-$.  K.~W\l odarczyk \cite{wlo87},  F.~Jafari  \cite{jaf93}  and M.~Abate \cite{ab98} obtained partial analogs
of the  Julia-Carath\'eodory Theorem for the polydisk.    Condition (\ref{A}) is replaced by the analogous
condition (\ref{zm}); it still implies the existence of nontangential limits and an inequality like  Julia's inequality (\ref{jul7}), but in dimensions greater than $1$ it does {\em not} imply the existence of an angular gradient (the natural analog of angular derivative: Definition \ref{4.4}). 

We show in this paper that, in the case of the bidisk $\D^2$, there is nevertheless a rich differential structure for a Schur class function $\ph$ at points $\tau$ satisfying (\ref{zm}):
\begin{equation}
\label{zm}
\liminf_{\la\to\tau} \frac{1-|\ph(\la)|}{1-\|\la\|} < \infty.
\end{equation}
Here $\| \la \|$ means $\max(|\la^1 |, | \la^2|)$.

 Any  $\ph$ satisfying (\ref{zm}), though not necessarily having an angular gradient at $\tau$, does have a directional derivative at $\tau$ in any direction pointing into the bidisk. 
 These directional derivatives need not vary linearly with the direction, as an angular gradient would,
but they do vary {\em  holomorphically} with the direction, and can be described by a function of one complex variable belonging to the Pick class (Theorem \ref{pickfn}).  Furthermore,  $\ph$ {\em can} have an angular gradient at $\tau$ even when not regular at that point.   When $\ph$ does have an angular gradient at $\tau$ we say that $\tau$ is a {\em $C$-point} of $\ph$.  The main result of the paper is that an analog of Carath\'eodory's result (C$\Rightarrow$D in Theorem \ref{thma1}) holds at $C$-points of functions in the Schur class of the bidisk.  That is, if $\tau\in\tb$ is a $C$-point of $\ph$ then the gradient of $\ph$ at $\la$ tends to the angular gradient of $\ph$ at $\tau$ as $\la$ tends nontangentially to $\tau$ (Theorem \ref{5.3p}).  A converse also holds.

To obtain these results we use Hilbert space models of  functions on  $\D^2$.  In this context a {\em model} of a function $\ph$ is an analytic map $u$ from $\D^2$ to a Hilbert space $\M$ such that a certain identity holds (see Definition \ref{1.4}).  Passage from a scalar-valued to a vector-valued function is not at first sight a gain in simplicity, but D. Sarason has shown in a beautiful monograph \cite{sar94} that the analogous process in one variable leads both to some simple proofs of classical results, including the Julia-Carath\'eodory Theorem, and to new insights.  Here we develop the model theory of functions at points for which  condition (\ref{zm}) holds; we call such points {\em $B$-points}.  On the way to proving the results described above we obtain new proofs of results of Jafari and Abate as they relate to the bidisk.

A limitation of models is that they are restricted to one and two dimensions, at least as far as the Schur class is concerned.
 For $n \geq 3$ it is not true that all functions in the Schur class possess models; those that do are said to belong to the {\em Schur-Agler class}.  There are many papers on interpolation by functions in the Schur-Agler class (for example \cite{agmc_bid, babo02}).  The present results too extend straightforwardly to the Schur-Agler class on the polydisk, but since there are already versions of Julia-Carath\'eodory for the whole Schur class of $\D^n$ it seems appropriate to confine ourselves to the case $n=2$ in this paper.

The paper is arranged as follows.  Section \ref{results} contains a statement of the Julia-Carath\'eodory theorem in one variable, makes precise two notions of nontangential convergence on $\D^2$, provides formal definitions of $B$-points and $C$-points  and presents the two  main results of the paper, Theorems \ref{pickfn} and \ref{5.3p}.  Section \ref{models} defines a model of a function and a realization of a model.  In Section \ref{JuliaLem} we state and prove the two-variable analogs of Julia's inequality and other results on nontangential convergence at a $B$-point.  In Section \ref{Bpoints} we characterize $B$-points in terms of models and prove the remarkable fact that if a model converges weakly along a nontangential sequence at a $B$-point then it also converges strongly.  We also introduce the important notion of the cluster set of a model at a $B$-point.
In Section \ref{example} we discuss the rational inner function
\beq
\label{eqa55}
\psi(\lambda) = \frac{ \frac{1}{2} \lambda^1 + \frac{1}{2}\lambda^2 -\la^1\la^2}{1 - \frac{1}{2}
\lambda^1 - \frac{1}{2}\lambda^2}
\eeq
which shows that Carath\'eodory's implication (A)$\Rightarrow$(C) of Theorem \ref{thma1} does not generalize to the bidisk.  Section \ref{directional} analyses the nature of directional derivatives at a $B$-point, showing how a function in the (one-variable) Pick class is associated with any $B$-point of a function in the Schur class of $\D^2$ (see Theorem \ref{pickfn}).  In Section \ref{holdiffC} we characterize $C$-points in terms of models and in Section \ref{carathm} we prove our Carath\'eodory theorem, Theorem \ref{5.3p}.
Section \ref{conclud} compares our results with those of Jafari and Abate.

\blue
We are very grateful to Zibin Huang for bringing to our attention an
error in the published version of this paper in {\em Mathematische Annalen} (2012).
In Propositions \ref{2.4} and \ref{Bptrealzn},
in Theorem \ref{4.19}, and in  Corollaries \ref{2.10} and \ref{4equiv}
 we incorrectly assumed that $\lim_{n \to
\infty} |\phi(\lambda_n)| = 1$ followed from the other hypotheses;
instead it needs to be assumed as a separate hypothesis.
We correct that mistake in this version. We also add a little more
explanation to the proof of Lemma 8.4.
New material is in blue.
\black

\section{Principal results} \label{results}
Here is the classical Julia-Carath\'eodory Theorem.
For $\tau\in\T$ the notation $\la \nt \tau$ means that $\la$ tends to $\tau$ along some set in $\D$ that approaches $\tau$ nontangentially; see Subsection \ref{nontang} below for a precise explanation.
\begin{theorem}
\label{thma1}
Let $\phi : \D \to \D$ be holomorphic and nonconstant. \black Let $\tau$ be a point on the unit circle
$\T$. The following conditions are equivalent:
\begin{enumerate}
\item[\rm (A)] there exists a sequence $\{ \ln \}$ in $\D$ tending to $\tau$ 
such that
\beq
\label{eqa2}
\frac{1 - | \phi(\ln)|}{1- |\ln|}
\eeq
is bounded;
\item[\rm(B)] for every sequence $\{ \ln \}$ tending to $\tau$ nontangentially, the quotient {\rm(\ref{eqa2})} is bounded;
\item[\rm (C)]  the nontangential limit
\[
\omega \df \lim_{\la \nt \tau} \ph(\la)
\]
and the angular derivative
\[
\ph'(\tau) \df \lim_{\la \nt \tau} \frac{\ph(\la)-\omega}{\la-\tau}
\]
exist;  
\item[\rm (D)] 
there exist $\omega \in\T$ and $\eta\in\C$ such that $\ph(\la) \to \omega$
and $\ph'(\la) \to \eta$ as $\la \to \tau$ nontangentially.
\end{enumerate}

When any of these conditions hold the nontangential limit
\[
\alpha \df \lim_{\la \nt \tau}
\frac{1-|\ph(\la)|}{1-|\la|}
\]
exists and is positive, and 
\beq \label{valeta2}
\eta= \ph'(\tau) = \al \bar\tau \omega.
\eeq
\black  Furthermore, for all $\la\in\D$,
\begin{equation}
\label{jul7}
\frac{|\ph(\la) - \omega |^2}
{1-|\ph(\la)|^2}
\ \leq \
\alpha \ \frac {|\la - \tau|^2} 
{1-|\la|^2} .
\end{equation}
\end{theorem}
G. Julia \cite{ju20} proved inequality (\ref{jul7}) under the hypothesis that $\ph$ has a Taylor expansion of the form $\ph(\la)=\omega +\eta(\la-\tau) + \zeta(\la-\tau)^2+o((\la-\tau)^2)$ (where $|\omega|=1$) valid for $\la\in\D$.  C. Carath\'eodory \cite{car29} identified the correct condition (A), and proved the remaining assertions.  Proofs may be found in \cite{car54,   sar94}, chapters I and VI respectively.

\subsection{Non-tangential Approach} \label{nontang}

If $S \subset \mathbb{D}$ and $\tau \in \D^-$, we say that $S$ {\em approaches $\tau$ nontangentially} if $\tau \in S^-$ and there exists a constant $c$ such that, for all $\lambda \in S$,
\[ 
| \tau - \lambda | \le c ( 1 - |\lambda|).
 \] 
With this terminology we have in mind primarily the case that $\tau\in\T$, but it is convenient to allow also the possibility that $\tau\in\D$.   In the latter case it follows from the definition that $S$  approaches $\tau$ nontangentially if and only if $\tau \in S^-\subset \D$.

We shall make use of similar terminology for the bidisk: if $S \subseteq \mathbb{D}^2$ and $\tau \in (\D^-)^2$, we say that $S$ {\em approaches $\tau$ nontangentially} if $\tau \in S^-$ and there is a constant $c$ such that 
\begin{equation} \label{1.1}
||\tau-\la|| \le c (1-||\lambda||) 
\end{equation} 
for all $\lambda \in S$.   Here and throughout the paper  $||\cdot||$ on $\C^2$ denotes the $\ell^\infty$ norm: for $\la\in\C^2$,
\[
|| \la || = \max\{|\la^1|,|\la^2|\}.
\]
We always use superscripts to denote the coordinates of points in $\C^2$.  We define the smallest $c>0$ for which the inequality (\ref{1.1}) holds to be the {\em aperture} of $S$. 

Note that $1-||\lambda||$ is the Euclidean distance between $\lambda$ and $\partial(\mathbb{D}^2)$, the topological boundary of $\mathbb{D}^2$,  so that the relation (\ref{1.1}) is the natural generalization to the bidisk of nontangential approach in the disk.

There is a second, more forgiving notion of nontangential approach. If $S \subseteq \mathbb{D}^2$, say that $S$ {\em approaches $\tau$ plurinontangentially} if $\tau  \in S^-$ and there exist sets $S^1, S^2 \subseteq \mathbb{D}$ such that $S^1$ approaches $\tau^1$ nontangentially, $S^2$ approaches $\tau^2$ nontangentially, and $S \subseteq S^1 \times S^2$. 

If $\{ \lambda_n\}$ is a sequence, we say that {\em $\lambda_n$ tends to $\tau$ nontangentially} or {\em plurinontangentially} if  the set
$ \{ \lambda_n : n \ge 1 \}$ approaches $\tau$ nontangentially or  plurinontangetially respectively. In these cases we write $\lambda_n \rightarrow \tau$ nt (or $\la_n \stackrel{\rm  nt }{\to} \tau$) and $\lambda_n \rightarrow \tau$ pnt.

\subsection{Two definitions and two theorems}\label{functions}

Let $\mathcal{S}$ denote the Schur class on the bidisk. Thus, $\phi \in \mathcal{S}$ means that $\phi$ is a holomorphic function on $\mathbb{D}^2$ and 
\[ 
||\phi||_\infty \df \sup_{\lambda \in \mathbb{D}^2} | \phi(\lambda)| \le 1. 
\]

\begin{definition} \label{2.1p}
Let $\phi \in \mathcal{S}$ and let $\tau \in  \partial (\D^2)$. \black We say that $\tau$ is
a {\em $B$-point} for $\phi$ if there exists a sequence $\{\lambda_n\}$ in $\mathbb{D}^2$ converging to $\tau$ such that 
\begin{equation} \label{2.3p}
\frac{ 1 - |\phi(\lambda_n)|}{ 1 -||\lambda_n||} \mbox{ is bounded. }
\end{equation}
\end{definition}

The following result was proved in \cite{jaf93,ab98}.  In our treatment it  follows from Corollaries~\ref{2.7} and \ref{2.10}. 
\begin{proposition}\label{thma2}
Let $\phi$ be in $\Schur$. The following are equivalent  for $\tau\in\partial (\D^2)$: 
\begin{enumerate}
\item[\rm (A)]  $\tau$  is a $B$-point for $\phi$;

\item[\rm (B)] for every sequence $\{\la_n\}$ in $\D^2$ that converges  nontangentially  to $\tau$
\[
\frac{1-|\ph(\ln)|}{1-\|\ln\|} \quad \mbox{ is bounded}.
\]
\end{enumerate}

When {\rm (A)} and {\rm (B)} are satisfied there exists $\omega\in\T$ such that
$\ph(\la) \to \omega$ as $\la \to \tau$ {\rm pnt}.
\end{proposition}

Of course every point of $\partial(\D^2)$ at which $\ph$ is regular is a $B$-point, but our concern is with $B$-points at which $\ph$ is {\em not} regular.  The example in Section \ref{example} shows that $\ph$ can even be singular in a {\em topological} sense at a $B$-point  $\tau$ (the function does not extend continuously to $\tau$).  In this case $\ph$ is not differentiable at $\tau$; nevertheless, $\ph$ {\em does}
 have a directional derivative at $\tau$ in all directions pointing into $\D^2$.  Our next theorem tells us that
the directional derivative can be expressed in terms of a function of a single complex variable.  

We shall denote by $\Pi$ the open upper half plane $\{z\in\C: \im z > 0\}$ and by $\mathbb{H}$ the open right half plane $\{z\in\C: \re z > 0\}$.
The {\em Pick class} is the class of analytic functions on $\Pi$ with non-negative imaginary part, that is
functions from $\Pi$ to $\Pi^-$.   It will be denoted by $\Pick$.
For $\tau\in\partial (\D^2)$ we define
\beq\label{defHtau}
\Htau = \left\{ \begin{array}{lcl}
    \tau^1\Ha\times \tau^2 \Ha & \mbox{ if } & \tau\in \T^2\\
     \tau^1\mathbb{H}\times \C & & \tau\in\T\times \D \\
    \C\times\tau^2\mathbb{H} & & \tau\in\D\times\T. \end{array} \right.
\eeq
For $\de\in\Htau$ the ray through $\tau$ in direction $-\de$ points into the bidisk and so it makes sense to consider the directional derivative
\[
D_{-\delta}\ph(\tau)= \lim_{t\to 0+} \frac{\ph(\tau-t\delta)-\ph(\tau)}{t}.
\]

Theorems~\ref{xdelta}, \ref{directderiv} and \ref{pickfunc} contain the following statement.
\begin{theorem}\label{pickfn}
Let $\tau \in\tb$ be a $B$-point of $\ph\in\Schur$.  For any $\de\in\Htau$ the directional derivative
$D_{-\de}\ph(\tau)$ exists, is an analytic function of $\de$ on $\Htau$ and is homogeneous of degree $1$ in $\de$.  
Moreover, if $\tau\in \T^2$, there exists a function $h$
such that both $h(z)$ and $-zh(z)$ are 
in the Pick class and analytic and real-valued on $(0,\infty)$,
\[
h(1) = -\liminf_{\la\to\tau}\frac{1-|\ph(\la)|}{1-\|\la\|}
\]
 and, for all $\de\in\Htau$, 
\beq \label{propg}
D_{-\de}\ph(\tau) = \ph(\tau)\overline{\tau^2}\de^2 h\left(\frac{\overline{\tau^2}\de^2}{\overline{\tau^1}\de^1}\right).
\eeq
\end{theorem}

An even stronger conclusion follows for $B$-points lying on $\tb\setminus\T^2$: roughly, $h(z)=-\mathrm{const}/z$ in equation (\ref{propg}) (Proposition \ref{Dphundist}).

In a forthcoming paper we plan to answer the inverse question: for which functions $h\in\Pick$ does there exist $\ph\in \Schur$ such that equation (\ref{propg}) holds?

Our second principal result applies to points in $\T^2$  at which $\ph$ has an {\em angular gradient} -- the natural analog for two variables of the angular derivative.
\begin{definition} \label{4.4}
Let $\phi \in \mathcal{S}, \tau \in \partial(\D^2)$.\\
{\rm(1)} For $S \subseteq \mathbb{D}^2$, $\tau \in S^-$ we say that $\phi$ {\em  has a holomorphic differential at $\tau$ on $S$} if there exist $\omega, \eta^1, \eta^2 \in \mathbb{C}$ such that, for all $\lambda \in S$,
\begin{equation} \label{4.5}
\phi(\lambda) = \omega + \eta^1(\lambda^1-\tau^1) + \eta^2(\lambda^2 - \tau^2) + e(\lambda) 
\end{equation}
where 
\begin{equation} \label{4.6}
\lim_{\lambda\rightarrow \tau, \, \, \lambda \in S} \frac{e(\lambda)}{||\lambda-\tau||} = 0.
\end{equation} 
{\rm(2)} We say that $\tau$ is  a {\em $C$-point for}  $\phi$  if, for every set $S$ that approaches $\tau$ nontangentially, $\ph$ has a holomorphic differential on $S$ and $\omega$ in relation {\rm(\ref{4.5})} is unimodular. \\
{\rm(3)} If $\tau$ is a $C$-point for $\ph$ we define the {\em angular gradient $\nabla\ph (\tau)$} of $\ph$ at $\tau$ to be the vector $\begin{pmatrix} \eta^1 \\ \eta^2\end{pmatrix}$, where $\ph$ has holomorphic differential {\rm (\ref{4.5})} on some set that approaches $\tau$ nontangentially.
\end{definition}
It is clear that, when $\tau$ is  a  $C$-point for $\ph$, the quantities $\omega, \eta^1,\eta^2$ in equation (\ref{4.5}) are the same for every nontangential approach region $S$, and so the definition of $\nabla\ph (\tau)$ in (3) is unambiguous.  In this situation we define $\ph(\tau)$ to be $\omega$.  

An apparent drawback of  the above definition of $C$-point is that a condition must hold for {\em every} set $S$ that  approaches $\tau$ nontangentially.  However, we shall see in Remark \ref{clarify} below that the condition need only be checked for a single suitable set $S$.

Every $C$-point is a $B$-point (Proposition~\ref{4.7}), and in one variable the converse also holds, by Theorem~\ref{thma1}. 
The function $\psi$ of equation (\ref{eqa55}) shows that, for functions of two variables, {\em not every $B$-point is a $C$-point}.  The two-variable analogs of conditions (C) and (D) from Theorem~\ref{thma1} are nevertheless equivalent.

\begin{theorem} \label{5.3p}
 If $ \tau \in \tb$ is a $C$-point for $\phi \in \Schur$ then 
\[
\lim_{\la \stackrel{\rm  nt }{\to} \tau} 
\nabla \phi(\lambda)\ =\ \nabla \phi(\tau).
\]
 Conversely, if $\tau\in\tb$ is a $B$-point of $\ph$ and $\lim_{\la \stackrel{\rm  nt }{\to} \tau} \nabla \phi(\lambda)$ exists then $\tau$ is a $C$-point of $\ph$. 
\end{theorem}
Points at which $\ph$ is regular are of course $C$-points.  The assertion of the theorem is trivial for such $C$-points,  but there are examples of functions in $\Schur$ that have singular $C$-points.  One example, to be constructed in a future paper, is the rational inner function
\[
\ph(\la)= \frac{-4\la^1(\la^2)^2+(\la^2)^2+3\la^1\la^2-\la^1+\la^2}{(\la^2)^2-\la^1\la^2-\la^1-3\la^2+4},
\]
which has a $C$-point at $(1,1)$, despite being singular there ($\ph$ cannot be extended continuously to $\D^2\cup \{(1,1)\}$).

The proofs of the theorems of this section rely on models and realizations of functions in the Schur class.  We discuss these in the next section.

\section{Models and Realizations} \label{models}

In one variable there is a model for any $\phi$ in the Schur class in terms of a Hilbert
space ${\mathcal H}(\ph)$, the {\em de Branges-Rovnyak} space for $\phi$, which
 is defined to be the Hilbert space of holomorphic
functions on  $\D$ with reproducing kernel
\[
k^\phi_{\la}(z) \ = \ \frac{1-\overline{\phi(\la)} \ph(z)}{1 - \bar \la z} \ = \
\langle k^\phi_\lambda, k^\phi_z \rangle_{{\mathcal H}(\phi)}
\]
It is convenient for us to introduce the function
$\psi(z) := \overline{\phi(\bar z)}$, and to use the vectors
\[
u_\la(z) \ = \ k^\psi_{\bar \la}(z) \ = \ \frac{1 - \overline{\psi(\bar \la)} \psi(z)}{1 - \la z} 
\]
which are analytic as a function of $\la$ (with values in ${\mathcal H}(\psi)$).
We then have the relation
\begin{equation}\label{eqb7}
1 - \overline{\phi(\mu)}\phi(\lambda) = (1 - \overline{\mu} \lambda) \langle  u_\lambda,u_\mu \rangle;
\end{equation}
we call the pair $({\mathcal H}(\psi), u)$ a {\em model} of $\ph$.

In two variables, even though no constructive formula is known for a model, there nevertheless does exist a model of any function in the Schur class, as was shown
by one of the authors \cite{ag90}:
for any $\phi \in \mathcal{S}$, there exists a separable Hilbert space $\mathcal{M}$, an orthogonal decomposition of $\mathcal{M}$, 
\[
\mathcal{M} = \mathcal{M}^1 \oplus \mathcal{M}^2, 
\] 
and an analytic map $u: \mathbb{D}^2 \rightarrow \mathcal{M}$ such that 
\begin{equation}\label{1.3}
1 - \overline{\phi(\mu)}\phi(\lambda) = (1 - \overline{\mu^1} \lambda^1) \langle  u_\lambda^1,u_\mu^1
\rangle  + (1 - \overline{\mu^2} \lambda^2)\langle u_\lambda^2,u_\mu^2 \rangle  
\end{equation}
 for all $\lambda, \mu \in \mathbb{D}^2$.
In equation (\ref{1.3}) we have written $u_\la$ for $u(\lambda)$, $u_\lambda^1 = P_{\mathcal{M}^1} u_\lambda$, and $u_\lambda^2 = P_{\mathcal{M}^2} u_\lambda$. 

In general, if $v \in \mathcal{M}$, we set $v^1 = P_{\mathcal{M}^1} v$ and $v^2 = P_{\mathcal{M}^2} v$. If $\lambda \in \mathbb{D}^2$, we may regard $\lambda$ as an operator on $\mathcal{M}$ by letting 
\[ 
\lambda v = \lambda^1 v^1 + \lambda^2 v^2
\] 
for $v \in \mathcal{M}$. 
 Note that $||\la||$ coincides with the operator norm of the operator
$\la$ acting on $\M_1 \oplus \M_2$.
With these notations, equation (\ref{1.3}) can be written in the slightly more compact form 
\[ 
1 - \overline{\phi(\mu)}\phi(\lambda) = \langle  (1 - \mu^* \lambda) u_\lambda, u_\mu \rangle . 
\]

We formalize the above notions into the following definition.
\begin{definition} \label{1.4}
Let $\phi \in \mathcal{S}$. We say that $(\mathcal{M}, u)$ is a {\em model of} $\phi$ if $\mathcal{M} = \mathcal{M}^1 \oplus \mathcal{M}^2$ is an orthogonally decomposed separable Hilbert space and $u: \mathbb{D}^2 \rightarrow \mathcal{M}$ is an 
analytic map such that equation {\rm (\ref{1.3})} holds for all $\lambda, \mu \in \mathbb{D}^2$. 
\end{definition}

Note that we can suppose without loss that $\{u_\la^j: \la\in\D^2\}$ spans a dense subspace of $\M_j$, since otherwise we may replace $\M_j$ by a subspace.  However it need not be the case that $\{u_\la:\la\in\D^2\}$ spans a dense subspace of $\M$.
\black

From the Hilbert space gadget $(\mathcal{M}, u)$ of Definition \ref{1.4} associated with a $\phi \in \mathcal{S}$ it is not difficult to go a step further and in turn attach to $(\mathcal{M}, u)$ an operator-theoretic construct. A lurking isometry argument yields the following result \cite{ag90}.
\begin{theorem} \label{1.5}
If $(\mathcal{M}, u)$ is a model of $\phi \in \mathcal{S}$, then there exist  $a\in\C$, vectors $\beta, \gamma \in \mathcal{M}$ and a linear operator $D: \mathcal{M} \rightarrow \mathcal{M}$ such that the operator
\beq\label{abcd}
\left[ \begin{array}{cc} a & 1\otimes\beta \\ \ga\otimes 1 & D \end{array} \right]
\eeq
is a contraction on $ \C \oplus \mathcal{M}$ and,  for all $ \la\in\D^2$,
\begin{eqnarray}
\label{eqaa6}
 (1-D \lambda) u_\lambda &=& \gamma,\\
   \ph(\la)  &=& a + \langle \la u_\la, \beta  \rangle .
\label{eqaa7}
\end{eqnarray}
 \end{theorem}

Any ordered 4-tuple $(a,\beta,\gamma,D)$ for which equations (\ref{eqaa6}) and
(\ref{eqaa7}) hold and for which the operator (\ref{abcd}) is a contraction will be called a {\em realization} of the model $(\M,u)$.

The concept of a realization of a model is more restricted than the
system-theoretic notion of a realization of a function.  Certainly, if equations (\ref{eqaa6}) and (\ref{eqaa7}) hold then
\beq \label{realiznform}
\ph(\la)= a + \langle \la(1-D\la)^{-1}\ga,\beta  \rangle ,
\eeq
 for all $\la\in\D^2$, and so $(a,\beta,\gamma,D)$ is a realization of $\ph$ in the sense of systems theory.  Consider, though, the (system-theoretic) realization $(0,
\beta,\gamma,0)$ of the function $\ph(\la)=\half\la^1$ where $\beta=[\half \quad
0]^T, \gamma=[1 \quad 0]^T$.  Equation (\ref{realiznform}) holds, and it is tempting to infer
that $(\C \oplus \C, u)$ is a model of $\ph$, where $u_\la=
(1-D\la)^{-1}\gamma =\gamma$.  However, it is simple to check that
equation (\ref{1.3}) does not hold.  Thus $(0,\beta,\gamma,0)$ is not a
realization of a model of $\ph$ in the sense of our definition.  To avoid
confusion we shall speak only of realizations of models, not functions.

\section{Julia's Lemma} \label{JuliaLem}
Julia's inequality (\ref{jul7}), the last part of Theorem \ref{thma1},  was extended to polydisks by K. W\l odarczyk \cite[Lemma 2.1]{wlo87}, F. Jafari \cite[Theorem 5]{jaf93}     and M. Abate \cite[Theorem 3.1]{ab98}.   W\l odarczyk obtained a version of inequality (\ref{jul7}) for the unit ball of any $J^*$-algebra, which includes the polydisk as a special case, while Jafari and Abate obtained analogs of the other parts of Theorem \ref{thma1}  for polydisks.  In this section, in the course of developing the theory of models, we give an alternative proof of Julia's inequality for the bidisk.

We define the {\em cluster set} of a model $(\mathcal{M},u)$ for a 
Schur class function at a $B$-point $\tau$ 
to be the set of limits in $\M$ of the weakly convergent sequences $\{u_{\la_n}\}$ as $\{\la_n\}$ ranges over all sequences in $\D^2$ that tend 
to $\tau$ and satisfy 
\beq \label{Cquo}
\frac{1-|\ph(\ln)|}{1-\|\ln\|} \mbox {  is bounded.  }
\eeq
  We shall denote the cluster set at $\tau$ by $Y_\tau$. \blue It follows from Proposition \ref{getx}(3) that the cluster set of a model of a Schur-class function at a B-point is nonempty. 
\black

The following observation will play an important role throughout the paper.
\begin{proposition} \label{getx}
Let $\tau\in\partial (\D^2)$ be a $B$-point for $\ph\in\Schur$ and let
$(\M,u)$ be a model of $\ph$.  \\
{\rm(1)}  If $x\in Y_\tau$ and $|\tau^j|< 1$ for $j=1$ or $2$ then $x^j=0$.\\
{\rm(2)} There exists $\omega\in\T$ such that for all $x\in Y_\tau$  and $\la\in\D^2$,
\begin{equation} \label{2.8}
1 - \overline{\omega} \phi(\lambda) = \sum_{|\tau^j|=1}(1 - \overline{\tau^j} \lambda^j) \langle  u_\lambda^j, x^j\rangle.
\end{equation} 
{\rm(3)} If
\beq\label{defalpha}
\liminf_{\la\in\D^2, \ \la\to\tau} \frac{1-|\ph(\la)|}{1-||\la||} = \al.
\eeq
then there exists $x\in Y_\tau$ such that $||x||^2 \leq \al$.
\end{proposition}
\begin{proof}
(3) Pick a sequence $\{\la_n\}$ in $\D^2$ converging to $\tau$ such that
\beq\label{chooselan}
 \lim_{n\to\infty}  \frac{1-|\ph (\la_n)|}{1- ||\la_n||} = \al;
\eeq
since $\tau$ is a $B$-point, $\al<\infty$. Passing to a subsequence if necessary, we can suppose that  $\ph(\la_n)$ converges to a point $ \omega\in\D^-$.  Inequality (\ref{chooselan}) implies that $\omega\in\T$.  By equation (\ref{1.3}) we have, since $1-|\la_n^j|^2 \geq 1- ||\la_n||^2$,
\begin{align}\label{ulanest}
\frac{1-|\ph(\la_n)|^2}{1-||\la_n||^2} &= \frac {(1-|\la_n^1|^2) ||u_{\la_n}^1||^2 + (1-|\la_n^2|^2) ||u_{\la_n}^2||^2}{1-||\la_n||^2} \nn \\
     & \geq  ||u_{\la_n}||^2.
\end{align}
Thus $||u_{\la_n}||$ is bounded.  By the compactness and metrizability of closed balls in $\M$ in the weak topology, we can take a further subsequence to arrange that $u_{\la_n}$ converges weakly to some $x\in\M$, where $ ||x||^2 \leq \al$.   Clearly $x\in Y_\tau$
and equation \blue (\ref{defalpha}) holds. \\ \black

(1) and (2) Consider $x\in Y_\tau$.  Pick a sequence $\{\la_n\}$ in $\D^2$ converging 
to $\tau$ such that $u_{\ln}\to x$ weakly and relation (\ref{Cquo}) holds.  By passing to a subsequence we can arrange that $\ph(\ln)\to \omega$ for some $\omega\in\D^-$.  Since (\ref{Cquo}) holds we have $|\omega|=1$.   On  letting $\mu = \lambda_n, \ n\to \infty$ in equation (\ref{1.3}) we deduce that, for $\la\in\D^2$, 
\begin{equation} \label{2.8a}
1 - \overline{\omega} \phi(\lambda) = (1 - \overline{\tau^1} \lambda^1) \langle  u_\lambda^1, x^1\rangle  + (1 - \overline{\tau^2} \lambda^2) \langle u_\lambda^2,x^2\rangle.
\end{equation} 
If $\tau\notin \T^2$, say $|\tau^i| < 1$, then on setting $\la=\ln$ in equation (\ref{2.8a}) and letting $n\to \infty$ we find that 
\beq \label{xi0}
x^i=0  \quad \mbox{  whenever  }  \quad |\tau^i| < 1.
\eeq
Thus statement (1) follows, and equation (\ref{2.8}) holds for all $\la\in\D^2$.  Suppose also $y\in Y_\tau$ -- say $u_{\mu_n}\to y$ where 
$\mu_n 
\to
\tau$ 
and $\ph(\mu_n) \to \zeta$.  Put $\la=\mu_n$ in equation (\ref{2.8}) and let $n\to \infty$ to obtain $1-\overline\omega \zeta =0$, that is, $\omega=\zeta$.  Thus $\omega$ is uniquely determined and independent of the choice of $x\in Y_\tau$, and so (2) holds. 
\end{proof}

Here is Julia's inequality for the bidisk (\cite[Lemma 2.1]{wlo87}, \cite[Theorem 5]{jaf93}, \cite[Theorem 3.1]{ab98}).
\begin{theorem} \label{Jineq}
Let $\ph\in\Schur, \ \tau\in  \partial(\D^2) \black$.  If there is a sequence $\{\la_n\} $ in $\D^2$ converging to $\tau$ such that 
\beq\label{alpha}
 \lim_{n\to\infty}  \frac{1-|\ph (\la_n)|}{1- ||\la_n||} = \al < \infty
\eeq
then there exists $\omega\in\T$ such that, for all $\la\in\D^2$,
\beq\label{JuIn}
\frac{|\ph(\la)-\omega|^2}{1-|\ph(\la)|^2} \leq \al \max_{|\tau^j|=1} \frac{|\la^j-\tau^j|^2}{1-|\la^j|^2}.
\eeq
If $\ph$ is not constant then $\al>0$.
\end{theorem}
\begin{proof}
We prove the result for $\tau\in\T^2$; obvious modifications yield the general case.
\black
Pick a model $(\M,u)$ of $\ph$.  By Proposition \ref{getx} there exist $\omega\in\T, \ x\in Y_\tau$ such that $||x||^2\leq \al$ and equation (\ref{2.8}) holds for all $\la\in\D^2$.

If $\al=0$ then $x=0$ and so, by equation (\ref{2.8}), $\ph$ is the constant function $\omega$.

Consider a fixed $\la\in\D^2$ and let
\[
R=\max_{j=1,2} \frac{|\la^j - \tau^j|^2}{1-|\la^j|^2}.
\]
 From equation (\ref{2.8}) we have
\[
| 1 -  \overline{\omega}\phi(\lambda) |  
\leq |1 - \overline{\tau^1} \lambda^1| \cdot ||x^1|| \cdot ||u_\lambda^1|| + |1 - \overline{\tau^2}\lambda^2|  \cdot||x^2|| \cdot ||u_\lambda^2||,
\]
and so, by the Cauchy-Schwarz inequality,
\begin{align*}
|\ph(\la) - \omega|^2  & \leq \{||x^1||^2+||x^2||^2\} \{ |\tau^1-\la^1|^2  \ || u_\la^1||^2 + |\tau^2-\la^2|^2  \ ||u_\la^2||^2 \} \\
   & \leq ||x||^2 R\{(1-|\la^1|^2) \ ||u_\la^1||^2 + (1-|\la^2|^2) \ ||u_\la^2||^2\}  \\
  &\leq \al R(1-|\ph(\la)|^2) \\
  &\leq \al R,
\end{align*}
which is the desired inequality (\ref{JuIn}).
\end{proof}
Theorem \ref{Jineq} can be reformulated in the terminology of horospheres and horocycles.
Recall that a {\em horocycle} in $\D$ is a set of the form $E(\tau, R)$ for some $\tau \in\T$ and $R>0$, where
\[
E(\tau, R) = \left\{\la\in\D: \frac{|\la-\tau|^2}{1-|\la|^2} < R\right\}.
\]
We shall denote by $D(z,r)$ the Euclidean disk in $\C$ with centre $z$ and radius $r>0$. $E(\tau,R)$ is the circular disk internally tangent to $\T$ at $\tau$ having radius $R/(R+1)$:
\[
 E(\tau,R) = D\left(\frac{\tau}{R+1}, \frac{R}{R+1}\right).
\]

For $\tau\in\D$ and any $R>0$ we define $E(\tau,R)$ to be $\D$.
\black
For $\tau\in \partial(\D^2)\black $ and $R>0$ the {\em horosphere} $E(\tau,R)$ is defined to be the set $E(\tau^1,R)\times E(\tau^2,R)$.  The following is then simply a restatement of inequality (\ref{JuIn}).
\begin{corollary} \label{Jineq2}
Under the hypotheses of Theorem {\rm \ref{Jineq}} there exists $\omega\in\T$ such that, for any $R>0$,
\beq\label{includ}
\ph(E(\tau,R)) \subset E(\omega,\al R).
\eeq
\end{corollary}

Another consequence of Theorem \ref{Jineq} is that one can test whether $\tau$ is a $B$-point using only the values of $\ph$ along the radius through $\tau$. 
\begin{corollary} \label{3.9}
Let $\ph\in\Schur$ and $\tau\in  \partial(\D^2)$.  The following conditions are equivalent.
\begin{enumerate}
\item[\rm (1)] $\tau$ is a $B$-point for $\ph$;
\item[\rm(2)] $(1-|\ph(\la)|)/(1-||\la||)$ is bounded on the radius $\{r\tau: 0<r<1\}$;
\item[\rm (3)] $(1-|\ph(\la)|)/(1-||\la||)$ is bounded on every set in $\D^2$ that approaches $\tau$ nontangentially.
\end{enumerate}  
  Moreover, if $\tau$ is a $B$-point for $\ph$,
\beq\label{radial}
 \lim_{r\to 1-}\frac{1-|\ph(r\tau)|}{1-r} = \liminf_{\la\to\tau, \, \la\in\D^2} \frac{1-|\ph(\la)|}{1-||\la||}.
\eeq
\end{corollary}
\begin{proof}
It is trivial that (3)$\Rightarrow$(2)$\Rightarrow$(1).   We prove (1)$\Rightarrow$(3).  Suppose $\tau$ is a $B$-point of $\ph$ and consider a set $S\subset \D^2$ that approaches $\tau$ with aperture $c\geq 1$.  Let $\al$ denote the lim inf on the right hand side of equation (\ref{radial}); by hypothesis, $\al$ is finite.   By Corollary \ref{Jineq2} there exists $\omega\in\T$ such that, for all $R> 0$,
\beq \label{Dineq}
\ph(E(\tau, R)) \subset E(\omega,\al R) = D\left( \frac{\omega}{\al R+1},  \frac{\al R}{\al R+1}\right).
\eeq
Pick any $\la\in S$ and $\ep > 0$ and let
\[
R = (1+\ep)\frac{\|\la-\tau\|^2}{1-\|\la\|^2}.
\]
Then $\la\in E(\tau, R)$ and, since $S$ has aperture $c$,
\[
0 < R \leq (1+\ep)\frac{c^2(1-\|\la\|)^2}{1-\|\la\|^2} \leq (1+\ep)c^2(1-\|\la\|).
\]
By the relation (\ref{Dineq})
\begin{align*}
 1- |\ph(\la)| &\leq |\ph(\la)-\omega| \\
  &\leq \frac{2\al R}{\al R+1}  \leq 2\al R \\
  &\leq 2(1+\ep) \al c^2(1-\|\la\|).
\end{align*}
Since this inequality holds for all $\ep > 0$ we have
\beq \label{caracondn}
\frac{1-|\ph(\la)|}{1-\|\la\|} \leq 2\al c^2.
\eeq
Hence (1)$\Rightarrow$(3).

We now prove equation (\ref{radial}).
By Theorem \ref{Jineq},  there exists $\omega\in\T$ such that the inequality (\ref{JuIn}) holds for all $\la\in\D^2$.  In particular, choosing $\la=r\tau$, we find that
\[
\frac{|\ph(r\tau)- \omega|^2}{1-|\ph(r\tau)|^2} \leq \al\frac{1-r}{1+r}\quad \mbox{ for } 0\leq r < 1.
\]
That is, $\ph(r\tau)$ lies in the horocycle $E(\omega,R)$ where $R= \al (1-r)/(1+r)$,  which is a circular disc of radius 
\[
\frac{R}{R+1} = \frac{\al (1-r)}{1+r+\al (1-r)}.
\]
Since $1-|\ph(r\tau)| \leq |\ph(r\tau)-\omega| \leq 2R/(R+1)$, we have
\[
\frac{1-|\ph(r\tau)|}{1-r} \leq \frac{|\ph(r\tau) -\omega|}{1-r} \leq  \frac{2\al }{1+r+\al (1-r)}.
\]
The right hand side is bounded for $0\leq r<1$, and has limit $\al$ as $r \to 1-$.  Thus
\[
\al = \liminf_{\la\to\tau, \, \la\in\D^2} \frac{1-|\ph(\la)|}{1-||\la||} \leq
\liminf_{r\to 1-}\frac{1-|\ph(r\tau)|}{1-r}\leq \limsup_{r\to 1-}\frac{1-|\ph(r\tau)|}{1-r} \leq \al .
\]
Equation (\ref{radial}) follows.
\end{proof}
Similar arguments to the above are given in \cite{jaf93, ab98}, Theorems 5 and 3.1 respectively.
 
A third consequence of Theorem \ref{Jineq} is that a Schur class function has a value at any $B$-point.  The following is contained in \cite[Theorem 3.1]{ab98}, where the terminology of ``restricted $E$-limit" is used. 
\begin{corollary}\label{2.7}
If $\phi \in \mathcal{S}$ and $\tau\in\partial (\D^2)$ is a $B$-point for $\phi$ then there exists $\omega \in \mathbb{T}$ such that $\phi(\lambda) \rightarrow \omega$ as $\lambda \rightarrow \tau$   horospherically, and {\em a fortiori} as $\lambda \rightarrow \tau$ plurinontangentially.
 \end{corollary} 

Here the {\em horospheric} topology on the closed bidisk is the topology for which a base consists of all open sets of $\D^2$ together with all sets of the form $\{\tau\}\cup E(\tau,R)$ where $\tau\in\partial (\D^2), \ R> 0$.    The horospheric topology is not Hausdorff on the closed bidisk (there do not exist disjoint neighborhoods of $\tau\in\T\times\D$ and any point of $\{\tau^1\}\times \D^-$), nor even $T_1$ (the closure of the singleton set $\{(1,0)\}$ is the face $\{1\}\times \D$), though it is Hausdorff on $\D^2\cup\T^2$.
\begin{proof}
Pick $\omega\in\T$ as in Corollary \ref{Jineq2}. 
We must show that, for every $\ep>0$, there exists a horospheric neighborhood $E(\tau, R)$ of $\tau$ such that 
\beq \label{ED}
\ph(E(\tau,R))\subset D(\omega, \ep).
\eeq
  In view of the inclusion (\ref{includ}) it suffices to choose $R'>0$ such that $E(\omega,R') \subset D(\omega, \ep)$ and then to take $R=R'/\al$.  Thus $\ph(\la)\to\omega$ as $\la\to\tau$ with respect to the horospheric topology.

Consider a net $(\ln)$ in $\D^2$ that tends to $\tau$ plurinontangentially, so that there is a set $S^j$ that approaches $\tau^j$  nontangentially such that all $\la_n^j \in S^j$ (or in the case that $\tau^j\in\D$, $S^j$ is relatively compact in $\D$) and $\ln\to \tau$.   We wish to show that $\ph(\ln)\to \omega$.

Let $\ep>0$ and choose  $R>0$ such that (\ref{ED}) holds.  We claim that there exists $\de >0$ such that, for $j=1,2$,
 \[
D(\tau^j,\de) \cap S^j \subset E(\tau^j, R).
\]
If $|\tau^j|<1$ then this inclusion is trivially satisfied, since $ E(\tau^j, R)=\D$.
Otherwise, if $S^1, S^2$ have apertures no greater than $c>1$, it suffices to take any $\de$ such that
\[
\de <  \min\left\{c, \frac{2Rc}{c^2+R} \right\}.
\]
Since $\ln\to\tau$ pnt, for large enough $n$ we have $\la_n^j\in D(\tau^j,\de) \cap S^j \subset E(\tau^j,R)$. Hence $\ln\in E(\tau,R)$ and so $\ph(\ln)\in D(\omega, \ep)$ as required.
\end{proof} \black

 With Corollary \ref {2.7} in mind we write simply $\phi(\tau) = \omega$  to mean (whenever $\phi \in \mathcal{S}$ and $\tau$ is a $B$-point for $\phi$) that there exists a sequence $\{\lambda_n\} \subseteq \mathbb{D}^2$ such that $\lambda_n \rightarrow \tau$ pnt and $\phi(\lambda_n) \rightarrow \omega$.

\section{$B$-points and models} \label{Bpoints}

In this section we will give a Hilbert space characterization of $B$-points in terms of models of $\phi$. 
This characterization leads naturally to a number of results about the interplay between the function theory and the model theory of $B$-points.  Recall Definition \ref{2.1p}:
 $\tau$ is a  $B$-point for $\phi$ if there exists a sequence $\{\lambda_n\}$ in $\mathbb{D}^2$  converging to $\tau$ such that 
\begin{equation} \label{2.3}
\frac{ 1 - |\phi(\lambda_n)|}{ 1 -||\lambda_n||} \mbox{ is bounded. }
\end{equation} 
\blue
\begin{newremark}[5.1a]\customlabel{5.1a}{limit-mod}
Let $\ph\in\mathcal{S}$, let $ \ \tau \in \partial(\D^2)$ and
let $\{ \lambda_n\}$ be a sequence in $\mathbb{D}^2$ such that $\lambda_n \rightarrow \tau$ nontangentially.  If condition \eqref{2.3} holds, then $\lim_{n\to\infty} |\ph(\la_n)| =1$.
For condition \eqref{2.3} implies that, for some $M\geq 0,$
\[
1-|\ph(\la_n)| \leq M(1-\|\la_n\|) \quad \text{for all} \quad n\in\mathbb{N}.
\]
Since $\la_n \to \tau \in \partial \D^2$, we have $\|\la_n\| \to 1$ and so $1-|\ph(\la_n)| \to 0$ as $n\to\infty$.

\end{newremark}

\black

The following proposition gives a criterion in terms of models for condition (\ref{2.3}) to hold.

\begin{proposition} \label{2.4}
Let $\phi \in \mathcal{S}.$
\blue
{\rm (i)} If $ \ \tau \in \partial(\D^2)\setminus\T^2$ \black
and $\{ \lambda_n\}$ is a sequence in $\mathbb{D}^2$ such that $\lambda_n \rightarrow \tau$ nontangentially then the following statements are equivalent. 
\begin{enumerate}
\item[\rm(1)] Condition {\rm (\ref{2.3})}  holds; 
\item[\rm(2)] \blue $ \lim_{n \to \infty} |\phi(\lambda_n)| =1$ and \black there exists  a  model $(\mathcal{M},u)$ of  $ \phi$   such  that $u_{\lambda_n}$ is  bounded; 
\item[\rm(3)] \blue $ \lim_{n \to \infty} |\phi(\lambda_n)| =1$ and \black for every model $(\mathcal{M}, u)$ of $ \phi, \, u_{\lambda_n}$  is  bounded. \black
\end{enumerate}
{\rm(ii)} If $ \ \tau \in \T^2$ 
and $\{ \lambda_n\}$ is a sequence in $\mathbb{D}^2$ such that $\lambda_n \rightarrow \tau$ nontangentially then the following statements are equivalent.
\begin{enumerate}
\item[\rm(1)] Condition {\rm (\ref{2.3})}  holds; 
\item[\rm(2)] there exists  a  model $(\mathcal{M},u)$ of  $ \phi$   such  that $u_{\lambda_n}$ is  bounded; 
\item[\rm(3)]  for every model $(\mathcal{M}, u)$ of $ \phi, \, u_{\lambda_n}$  is  bounded.\black
\end{enumerate}

\end{proposition}
\begin{proof} 
Fix $c$ such that $||\tau - \lambda_n|| \le c  (1 - ||\lambda_n||)$ for all $n$.

\blue
Clearly, in both cases (i) and (ii), (3) implies (2).  

Let us show that in both cases (i) and (ii), (1)$\Rightarrow$(3). Assume (1): say condition (\ref{2.3}) holds with bound $M$. 
\black
If $(\mathcal{M}, u)$ is a model of $\phi$, then for $\la=\ln$
\begin{align*}
( 1-||\la||)  \ ||u_{\la}||^2 & \le  ( 1 - |\la^1|) \ ||u^1_{\la}||^2 + (1-|\la^2| ) \ ||u^2_{\la}||^2 \\
& \le (1 - |\la^1|^2) \ || u_{\la}^1||^2 + (1 - |\la^2|^2) \ ||u_{\la}^2||^2 \\
& = 1 - |\phi(\la)|^2 \\
& \le  2 M( 1 - ||\la||).
\end{align*} \blue
Hence, in both Cases,  (1) implies that, for all models $(\mathcal{M},u)$ of  $ \phi$, $u_{\lambda_n}$ is  bounded. By Remark \ref{limit-mod},   $ \lim_{n \to \infty} |\phi(\lambda_n)| =1$. Therefore  (1)$\Rightarrow$(3). \black

It remains to show that (2)$\Rightarrow$ (1).
In Case (ii), assume (2):  there exists  a  model $(\mathcal{M},u)$ of  $ \phi$   such  that $u_{\lambda_n}$ is  bounded,
\black  say $||u_{\la_n}|| \leq M$ for some model of $\ph$.  
Suppose that $\tau\in\T^2$.  For $\la=\la_n$ we have
\begin{align} \label{star} 
1 - |\phi(\la)|   & \le  1 - | \phi(\la)|^2   \nn\\
& = (1 - |\la^1|^2)  || u^1_{\la}||^2 \ + \  (1-|\la^2|^2)  ||u_{\la}^2||^2 \nn \\
& \le  2 (1 - |\la^1|)   || u^1_{\la}||^2 \ + \  2 (1-|\la^2|) ||u_{\la}^2||^2 \\
&  \le 2 | \tau^1 - \la^1| \cdot ||u^1_{\la}||^2 \ + \  2 |\tau^2 - \la^2 |\cdot ||u_{\la}^2||^2 \nn \\
& \le 2 c  ( 1 - ||\la||)||u_{\la}||^2  \leq 2cM^2 \, (1-||\la||). \nn
\end{align}
Thus (2)$\Rightarrow$(1) \blue in the case when $\tau\in\T^2$. \black

\blue In Case (i) consider $\tau\in \T\times \D$ (the argument for $\tau \in \D\times\T$ is similar) and assume (2) holds:
 $ \lim_{n \to \infty} |\phi(\lambda_n)| =1$ and  there exists  a  model $(\mathcal{M},u)$ of  $ \phi$   such  that $u_{\lambda_n}$ is  bounded.

Note that $ \lim_{n \to \infty} |\phi(\lambda_n)| =1$, and so there exist a subsequence of $\{\la_n\}$ and $\omega\in\T$  such that 
\[
 \lim_{n\to\infty} \ph(\la_n) =\omega.
\]
On  letting $\mu = \lambda_n, \ n=1,2 \dots,$ in equation (\ref{1.3}) we deduce that, for $\la\in\D^2$, 
\begin{equation*}\tag{5.3a} \label{2.8b}
1 - \overline{\phi(\lambda_n)} \phi(\lambda) = (1 - \overline{{\lambda_n}^1} \lambda^1) \langle  u_{\lambda}^1, 
u_{\lambda_n}^1\rangle  + (1 - \overline{{\lambda_n}^2} \lambda^2) \langle u_\lambda^2,u_{\lambda_n}^2\rangle.
\end{equation*} 
Since  $u_{\lambda_n}$ is  bounded in $\M$, we may pass to a subsequence of $\{\lambda_n\}$ for which there exists $x\in Y_\tau$ such that $u_{\lambda_n}$ converges weakly to the point $x=(x^1, x^2)$.
Take limits in equation \eqref{2.8b} as $\lambda_n \rightarrow \tau$, to obtain, for $\la\in\D^2$,
\begin{equation*}\tag{5.3b} \label{2.8c}
1 - \overline{\omega} \phi(\lambda) = (1 - \overline{\tau^1} \lambda^1) \langle  u_\lambda^1, x^1 \rangle  + (1 - \overline{\tau^2} \lambda^2) \langle u_\lambda^2,x^2 \rangle.
\end{equation*} 

Since $\tau \in \T\times\D$, $|\tau^2| < 1$, then on setting $\la=\la_n$ in equation (\ref{2.8c}) and letting $n\to \infty$ we find that $x^2=0$. Thus, by equation \eqref{2.8c}, for all $\la\in\D^2$, 
\begin{equation*} \tag{5.3c}\label{omegaphi}
\bar\omega \ph(\la) = 1+ \overline{\tau^1}(\la^1-\tau^1)\lan u_\la^1,x^1 \ran.
\end{equation*} \black
Hence, by equation \eqref{omegaphi},
\[
1-|\ph(\la)|^2 = \left| -2\re \overline{\tau^1}(\la^1-\tau^1)\lan u_\la^1,x^1 \ran - |\la^1-\tau^1|^2 |\lan u_\la^1,x^1 \ran|^2 \right|.
\]
For $\la$ close to $\tau$ we have $||\la||= |\la^1|$.  In addition,  for $\la = \la_n$ and large $n$, \ $|\la^1-\tau^1| < 1$ and so
\begin{align*}
\frac{1-|\ph(\la)|^2}{1-||\la||^2} &\leq  \frac {2|\la^1-\tau^1| \  |\lan u_\la^1,x^1 \ran | + |\la^1-\tau^1|^2  \ |\lan u_\la^1,x^1 \ran |^2}{1-|\la^1|} \\
  & \leq 2cM||x^1|| + cM^2 \ ||x^1||^2.
\end{align*}
Hence (2)$\Rightarrow$(1).
 \end{proof}

\blue
\begin{newremark}[5.3d]\customlabel{5.3d}{5.3d}
\rm
Note the difference between parts (i) and (ii) of Proposition \ref{2.4}.  Whereas for $\tau\in\T^2$, $\tau$ is a B-point for $\ph$ if and only if $u_{\lambda_n}$ is bounded for some model $(\M,u)$ of $\ph$ and every sequence $(\lambda_n)$ in $\D^2$ that converges nontangentially to $\tau$, for $\tau\in \partial(\D^2)\setminus \T^2$ this equivalence does not hold.  For consider the example $\ph(\lambda) = \lambda^2$ for $\lambda=(\lambda^1,\lambda^2) \in \D^2$.
A model for $\ph$ is $(\M,u)$ where $\M={0}\oplus\C$ and $u:\D^2 \to \M$ is given by $u(\la)=(u^1(\la),u^2(\la))$ where $u^1(\la)=0$ and $u^2(\la)=1$ for all $\la\in\D^2$. The model relation becomes
\[
1-\overline{\ph(\mu)}\ph(\la) = (1-\overline{\mu^1} \la^1)0+ (1-\overline{\mu^2} \la^2) \langle{1},{1} \rangle_\C,
\]
which is clearly valid.  Here $u(\la)$ is bounded, and the point $\tau=(1,0)$ is not a B-point of $\ph$, for if $\la_n$ is any sequence in $\D^2$ that converges to $\tau$ we have $\la_n^1 \to 1, \la_n^2 \to 0$ and
\[
\frac{1-|\ph(\la_n)|}{1-\|\la_n\|} = \frac{1-|\la_n^2|}{1-\max(|\la_n^1|,|\la_n^2|)},
\]
which is unbounded, since the numerator tends to $1$ and the denominator tends to $0$.

\end{newremark}

\black
 \begin{proposition} \label{2.9}
Let $\phi \in \mathcal{S}$ and let $(\mathcal{M},u)$ be a model of $\phi$. If $\tau\in\partial(\D^2)$ is a $B$-point for $\phi$ then $u_\lambda$ is bounded on any set in $\mathbb{D}^2$ that approaches $\tau$ nontangentially.   
In fact, if $S$ approaches $\tau$ with aperture $c>0$ and $\al$ is defined by equation {\rm (\ref{defalpha})} then
\beq \label{boundula}
||u_\la|| \leq 2 c \sqrt{\al} \quad \mbox{ for all  } \la\in S.
\eeq
\end{proposition}
 \begin{proof} By Corollary \ref{2.7} there exists $\omega \in \mathbb{T}$ such that $\phi(\lambda) \rightarrow \omega$ as
$\lambda \rightarrow \tau$ pnt.  By Proposition \ref{getx}, there is an $x \in \mathcal{M}$ such that $||x||^2 \leq \al$ and 
equation (\ref{2.8}) holds. 
Fix a set $S$ that approaches $\tau$ with aperture $c>0$. For $\lambda \in S,$
\begin{align*} 
(1 - ||\lambda||^2) ||u_\lambda||^2 
& \le  (1-|\lambda^1|^2)||u_\lambda^1||^2 + (1-|\lambda^2|^2)||u_\lambda^2||^2 \\
& = 1 - |\phi(\lambda)|^2 \\
& = 1 -\overline{\omega}\phi(\lambda) + \overline{(\omega - \phi(\lambda))} \phi(\lambda) \\
& \le 2 |1 -\overline{\omega}\phi(\lambda)| \\
& = 2 \sum_{|\tau^j|=1}| (1 - \overline{\tau^j}\lambda^j) \langle u_\lambda^j,x^j\rangle  | \\
& \le 2   \sum_{|\tau^j|=1} |\tau^j-\lambda^j| \ ||u_\lambda^j|| \ ||x^j||  \\
& \le 2c (1-||\lambda||) \ \sum_{|\tau^j|=1} ||u_\lambda^j||\  ||x^j|| \\
& \le 2c (1-||\lambda||)  \ ||u_\lambda|| \  ||x|| \\
& \le 2c (1 - ||\lambda||^2)\ ||u_\lambda|| \sqrt{\al}.
\end{align*} 
Hence, for $\lambda \in S, \ 
|| u_\lambda || \le 2c \sqrt{\al}$.
\end{proof}
 \begin{remark}\label{ulaface}
\rm  For a $B$-point $\tau\in\T\times\D$ the same argument gives a stronger boundedness property: for any $\si\in\{\tau^1\}\times \D$, $u_\la$ is bounded on any set $S$ that approaches $\si$ nontangentially.  The above reasoning yields $||u_\la|| \leq 2c||x^1||$ for $\la\in S$ when $S$ has aperture $c$.
\end{remark}
\black
Proposition \ref{2.9} has the following corollary, which sheds light on the nature of Definition \ref{2.1p}. 

 \begin{corollary} \label{2.10}
Let $\phi \in \mathcal{S}$, let $\tau \in \tb$ and let $(\mathcal{M},u)$ be a model of $\phi$.\\
\blue
{\rm (i)} If $ \ \tau \in \partial(\D^2)\setminus\T^2$, then  the following conditions are equivalent: \black
 \begin{enumerate}
\item [\rm(1)] $\ds \frac{1 - |\phi(\lambda_n)|}{1-||\lambda_n||}$ is bounded  on  some  sequence $\{\lambda_n\}$ that converges to $ \tau$;
\item[\rm(2)] $ u_{\lambda_n}$ is bounded  on  some  sequence $\{\lambda_n\}$  that  approaches $\tau$ nontangentially, and
\blue  $ \lim_{n \to \infty} |\phi(\lambda_n)| =1$; \black 

\item[\rm (3)] $\ds \frac{1 - |\phi(\lambda)|}{1-||\lambda||}$ is  bounded  on  every  subset  of $ \mathbb{D}^2$ that  approaches $\tau $ nontangentially;

\item[\rm(4)] for every subset $S$ of $\D^2$ that approaches $\tau$ nontangentially, \blue $\lim_{\la\to\tau, \la\in S} |\phi(\lambda)| =1$, and \black
$ u_\lambda $ is  bounded  on $S$.

\end{enumerate}
{\rm (ii)} If $ \ \tau \in \T^2$, then  the following conditions are equivalent: \black
 \begin{enumerate}
\item [\rm(1)] $\ds \frac{1 - |\phi(\lambda_n)|}{1-||\lambda_n||}$ is bounded  on  some  sequence $\{\lambda_n\}$ that converges to $ \tau$;
\item[\rm(2)] $ u_{\lambda_n}$ is bounded  on  some  sequence $\{\lambda_n\}$  that  approaches $\tau$ nontangentially;
\item[\rm (3)] $\ds \frac{1 - |\phi(\lambda)|}{1-||\lambda||}$ is  bounded  on  every  subset  of $ \mathbb{D}^2$ that  approaches $\tau $ nontangentially;
\item[\rm(4)] $ u_\lambda $ is  bounded  on  every  subset  of $\mathbb{D}^2 $ that  approaches $ \tau $ nontangentially.
\end{enumerate}
\end{corollary}
\begin{proof}  (1)$\Leftrightarrow$(2) and (3)$\Leftrightarrow$(4) are equivalent by Propositions \ref{2.4} and \ref{3.9}.   (3)$\Rightarrow$(1) and  (4)$\Rightarrow$(2) are trivial. Finally, (1) is just the definition that $\tau$ be a $B$-point for $\phi$, so that the implication (1)$\Rightarrow$(4) is just a restatement of Proposition \ref{2.9}. 
\end{proof}

Note that, since in the hypotheses of the corollary we may take $(\mathcal{M}, u)$ to be an arbitrary model of $\phi$, the corollary implies that (2) and (4) hold for some model of $\phi$ if and only if they hold for every model of $\phi$. 

In one variable, if $\tau$ is a $B$-point for a $\phi \in \mathcal{S}$ and $(\mathcal{M},u)$ is a model of $\phi$, then $u_\lambda$ extends to the point $\tau$ in such a way that  $u$ is continuous on $S\cup\{\tau\}$ for any set $S$ that approaches $\tau$ nontangentially; this follows from the approach of Sarason \cite{sar94}. We present an example in the next section that shows that this continuity phenomenon is absent in two variables. Nevertheless, the singularity of $u_\lambda$ at $\tau$ is quite tame, a fact which we explore for the remainder of this section. 
\begin{proposition} \label{2.15}
Let $\phi \in \mathcal{S}$, let $(\mathcal{M},u)$ be a model of $\phi$ and let $\tau\in \tb$ be a $B$-point for $\phi$. Suppose that $\{\lambda_n\}$ converges to $\tau$ nontangentially in $\mathbb{D}^2$.  If $\{u_{\lambda_n}\}$ converges weakly in $\M$ then $\{u_{\lambda_n}\}$ converges in norm. 
\end{proposition}
\begin{proof}
Fix a sequence $\{\lambda_n\}$ in $\mathbb{D}^2$ such that $\lambda_n \rightarrow \tau$ nt and
 $u_{\lambda_n} \rightarrow x$ weakly in $\mathcal{M}$.   Let $\omega=\ph(\tau)$.
By Proposition \ref{getx}, equation (\ref{2.8}) holds for all $\lambda \in \mathbb{D}^2$.   

As
\begin{equation}
\| u_{\ln} - x \|^2 \ = \ \left( \| u_\ln \|^2 - \Re \lan u_\ln , x \ran \right) +
\left( \|x \|^2 - \Re \lan u_\ln, x \ran \right),
\label{eqbf}
\end{equation}
and the second term in (\ref{eqbf}) tends to zero, it suffices to prove that
$ \| u_\ln \|^2 - \Re \lan u_\ln , x \ran $ tends to zero.

For any $\la \in \D^2$, Proposition~\ref{getx} tells us that
\begin{equation}
\label{ref1}
1 - \overline{\omega} \phi(\lambda) = \sum_{j=1}^2
(1 - \overline{\tau^j} \lambda^j) \langle  u_\lambda^j, x^j\rangle.
\end{equation}
So subtracting and adding twice the real part of (\ref{ref1}) 
(and using the defining property (\ref{1.3}) of a model)
we get the
first equality
in the following string; the rest follow from rearranging the terms.
\begin{eqnarray}
\lefteqn{ 
\nonumber
\sum_{j=1}^2 ( 1 - |\la^j|^2)  \left( \| u_\la^j \|^2 - \Re \lan u_\la^j , x^j \ran \right)}
\\
\nonumber
&&= \ 1 - |\phi(\la) |^2 - 2 \Re ( 1 - \bar\omega \phi(\la) ) 
+\
\Re \sum_{j=1}^2 
 [2(1 - \bar\tau^j \la^j) - (1 - |\la^j|^2)] \lan u_\la^j, x^j \ran \\
&&
\nonumber
= \
- |1-\bar \omega \phi(\la)|^2  +
\Re \sum_{j=1}^2 (1 - 2 \bar \tau^j \la^j  + |\la^j|^2 ) \lan u_\la^j, x^j \ran \\
&&
\nonumber
= \
- |1-\bar \omega \phi(\la)|^2 +
\Re \sum_{j=1}^2 (1 - 2 \bar \tau^j \la^j  + |\la^j|^2 ) \| x^j \|^2
\\
\label{longq}
 && \mbox{}\qquad 
+\  \Re \sum_{j=1}^2 (1 - 2 \bar \tau^j \la^j  + |\la^j|^2 ) \lan u_\la^j -x^j , x^j \ran
\end{eqnarray}
 We finish by considering two cases. 
First, assume $\tau$ is in $\T^2$. Then as $\ln \to \tau$ nt, the four quantities 
$$
1-|\la_n^1|^2, \
1-|\la_n^2|^2, \
1 - \| \ln \|, \
\| \tau - \ln \|
$$
are all comparable.
Therefore, dividing through in (\ref{longq}) gives
\begin{eqnarray*}
\lefteqn{
\| u_\ln \|^2 - \Re \lan u_\ln, x \ran }
\\&& \leq \
c \left[ \frac{|1 - \bar \omega \phi(\ln) |^2}{1 - \| \ln \|}
+ 
\sum_{j=1}^2 \frac{ \Re (1 - 2 \bar \tau^j \la_n^j  + |\la_n^j|^2 ) }{1 - |\la_n^j|^2} \| x \|^2 \right.
\\
&&\qquad + \left. \sum_{j=1}^2 \frac{ | 1 - 2 \bar \tau^j \la_n^j  + |\la_n^j|^2 | }{1 - |\la_n^j|^2}
| \lan u_\ln -x , x \ran | \right] \\
\\&& \leq \
c \left[ \frac{|1 - \bar \omega \phi(\ln) |^2}{1 - \| \ln \|}
+ 
\sum_{j=1}^2 \frac{ | \tau^j -\la_n^j |^2  }{1 - |\la_n^j|^2} \| x \|^2 \right.
\\
&&\qquad + \left. \sum_{j=1}^2 \frac{ |  \tau^j -\la_n^j  |(1+|\la_n^j|) }{1 - |\la_n^j|^2}
| \lan u_\ln -x , x \ran | \right].
\end{eqnarray*}
The first term on the right tends to zero by Corollary~\ref{2.10}, 
 the second because $\ln$ tends to $\tau$ nt, and in the third term each summand is the product of a bounded factor and a factor that tends to zero.

\vskip 5pt
Case: $\tau \in \T \times \D$.
By Proposition  \ref{getx}, we have $x^2 = 0$. Indeed, as
\[
1 - | \phi (\ln )|^2 \ = \ (1- | \la_n^1|^2) \| u^1_\ln \|^2 \ + \ 
(1- | \la_n^2|^2) \| u^1_\ln \|^2 ,
\]
letting $\ln$ tend to $\tau$ we get $\| u^2_\ln \|$ tends to zero.

Now repeat the calculations in (\ref{longq}), and notice that the left-hand side is less than or equal to the
term with just $j=1$.
Therefore
\begin{eqnarray*}
\lefteqn{
\| u^1_\ln \|^2 - \Re \lan u^1_\ln, x^1 \ran }
\\&& \leq \
c \left[ \frac{|1 - \bar \omega \phi(\ln) |^2}{1 - \| \ln \|}
+ 
 \frac{ \Re (1 - 2 \bar \tau^1 \la_n^1  + |\la_n^1|^2 ) }{1 - |\la_n^1|^2} \| x^1 \|^2 \right.
\\
&&\qquad + \left. \frac{ | 1 - 2 \bar \tau^1 \la_n^1  + |\la_n^1|^2 | }{1 - |\la_n^1|^2}
| \lan u^1_\ln -x^1 , x^1 \ran | \right],
\end{eqnarray*}
and again all the terms on the right tend to zero.
\end{proof}

\begin{theorem} \label{2.16}
Let $\phi \in \mathcal{S}$, let $\tau \in \partial (\D^2)$ be a $B$-point for $\phi$ and let $(\mathcal{M},u)$ be a model of $\phi$.   For any realization $(a,\beta,\ga,D)$ of $(\mathcal{M},u)$ there exists a unique vector $u_\tau \in \mathcal{M}$ such that $(1-D\tau)u_\tau = \gamma$  and $u_\tau \perp \ker(1-D\tau)$. Furthermore, if $S\subset \mathbb{D}$ approaches $1$ nontangentially, then 
\begin{equation}\label{2.17}
\lim_{z\rightarrow  1,  \  z \in S} u_{z\tau} = u_\tau. 
\end{equation}  
Consequently
\beq \label{normutau}
\|u_\tau\|^2 = \liminf_{\la\to \tau} \frac{1-|\ph(\la)|}{1-\|\la\|}.
\eeq 
\end{theorem}
\begin{proof} 
Choose a sequence $z_n \rightarrow 1$ such that $u_{z_n\tau} \rightarrow x \in \mathcal{M}$. From (\ref{eqaa6}) we conclude that $(1-D\tau)x = \gamma$. Since $D\tau$ is a contraction, $\ker(1-D\tau) \perp \mathrm{ran}(1-D\tau)$. Hence $\gamma \perp \ker(1-D\tau)$, and it follows from (\ref{eqaa6}) that,
for all $z \in \mathbb{D}$,
\[
u_{z\tau}= (1-D\tau)(1-zD\tau)^{-1}x,
\]
and so $u_{z\tau} \perp \ker(1-D\tau)$.
Hence $x \perp \ker(1-D\tau)$. We have shown that there is a vector $x \in \mathcal{M}$ with the properties that $(1-D\tau)x = \gamma$ and $x \perp \ker(1-D\tau)$. Since such a vector is unique we deduce the first assertion of the theorem by taking $u_\tau = x$. 

To see (\ref{2.17}), suppose that $z_n \rightarrow 1$ and $u_{z_n \tau} \rightarrow v$. Then $v \perp \ker(1-D\tau)$ and $(1-D\tau)v=\gamma$. Hence $v=u_\tau$. 

For any $r\in(0,1)$ we have, by (\ref{1.3}),
\[
\frac{1-|\ph(r\tau)|^2}{1-r^2} = \|u_{r\tau}\|^2.
\]
As $r\to 1-$ the right hand side tends to $\|u_\tau\|^2$, by equation (\ref{2.17}), and the left hand side tends to the lim inf in equation (\ref{normutau}), by Corollary \ref{3.9} and Proposition \ref{2.15}.  This establishes equation (\ref{normutau}). \black
\end{proof}

\begin{proposition} \label{2.18}
Let $\phi \in \mathcal{S}$, let $\tau \in \partial (\D^2)$ be a $B$-point for $\phi$ and let $\phi(\tau) =\omega$. If $(\mathcal{M},u)$ is a model of $\phi$ and $x\in Y_\tau$, the cluster set  of the model at $\tau$, then, for any realization $(a,\beta,\ga,D)$ of $(\M,u)$, there exists $e\in\ker(1 -D\tau)$ such that
\begin{equation} \label{2.20}
x= u_\tau +e
\end{equation}
 where $u_\tau$ is the vector described in Theorem {\rm\ref{2.16}}.
\end{proposition} 
\begin{proof} 
 By equation (\ref{eqaa6}), $(1-D\tau)x = \gamma$. Hence $x-u_\tau \in \ker(1-D\tau)$. 
\end{proof}

There is a simple characterization of $B$-points in terms of realizations.

\begin{proposition}\label{Bptrealzn}
Let $(\M,u)$ be a model of $\ph\in\Schur$ and let $(a, \beta,\gamma,D)$ be a realization of $(\M,u)$.  \\
{\rm (i)} The following conditions are equivalent for a point $ \ \tau \in \partial(\D^2)\setminus\T^2$:
\begin{enumerate}
\item [\rm(1)] $\tau$ is a $B$-point for $\ph$;
\item[\rm(2)] $\ga\in\mathrm{ran}~(1-D\tau)$ and \blue 
$ \lim_{\la \nt \tau} |\phi(\lambda)| =1$. \black 
\end{enumerate}
{\rm (ii)} The following conditions are equivalent for a point  $\tau\in\T^2$:
\begin{enumerate}
\item [\rm(1)] $\tau$ is a $B$-point for $\ph$;
\item[\rm(2)] $\ga\in\mathrm{ran}~(1-D\tau)$.
\end{enumerate}
\end{proposition}

\begin{proof} In both cases (1) implies (2), since by Theorem \ref{2.16},  if $\tau$ is a $B$-point then there exists $u_\tau\in\M$ such that $(1-D\tau)u_\tau=\ga$, and hence $\ga\in\mathrm{ran}~(1-D\tau)$. 
\blue
By Remark \ref{limit-mod}, if $\tau$ is a  $B$-point for $\phi$
 then $ \lim_{\la \nt \tau} |\phi(\lambda)| =1$.
\black

To prove that (2) implies (1), suppose that $(1-D\tau)x=\ga$ for some $x\in\M$.  Then, for $\la\in\D^2$,
\begin{eqnarray*}
u_\la &=&(1-D\la)^{-1}\ga =(1-D\la)^{-1}(1-D\tau)x \\
   &=&(1-D\la)^{-1}(1-D\la +D(\la-\tau))x\\ & =& x+(1-D\la)^{-1}D(\la-\tau)x.
\end{eqnarray*}
Since $D$ is a contraction, 
\[
|| (1-D\lambda)^{-1}D||  = || \sum_{n=0}^\infty (D\lambda)^n D|| 
  \le \sum_{n=0}^\infty ||\lambda||^n  = \frac{1}{1-||\lambda||},
\]
and so
\begin{eqnarray*}
||u_\la||  &\leq& ||x|| + ||(1-D\la)^{-1}D||\cdot||(\la-\tau)x||\\
  &\leq& ||x|| +\frac{||\la-\tau||}{1-||\la||} ||x||.
\end{eqnarray*}
Hence, for all $\la$ in a set $S$ that approaches $\tau$ nontangentially with aperture $c$,
\[
||u_\la|| \leq (1+c)||x||.
\]
Thus $u_\la$ is bounded on $S$. 
 \blue In the case when  $ \ \tau \in \partial(\D^2)\setminus\T^2$,   by assumption, 
 $ \lim_{\la \nt \tau} |\phi(\lambda)| =1$.  \black 
By Corollary \ref{2.10}, $\tau$ is a $B$-point of $\ph$.
\end{proof}

 \section{An Example} \label{example}
 In this section we present a simple example of a rational inner function on the bidisk that has a $B$-point at $(1,1) \in \mathbb{T}^2$ but does not satisfy the 2-dimensional analog of the conclusion
 $(B) \Rightarrow (C)$ of the \cjw theorem.

Let $\phi$ be the inner function
\beq
\label{eqd8}
\phi(\lambda) = \frac{ \frac{1}{2} \lambda^1 + \frac{1}{2}\lambda^2 -\la^1\la^2}{1 - \frac{1}{2}\lambda^1 - \frac{1}{2}\lambda^2}.
\eeq
 
A straightforward calculation yields 
\[
1 - \overline{\phi(\mu)} \phi(\lambda) = (1 - \overline{\mu^1}\lambda^1)\langle u_\lambda^1,u_\mu^1\rangle  +(1 - \overline{\mu^2}\lambda^2)\langle u_\lambda^2,u_\mu^2\rangle 
\]
where
\[
u_\lambda =  \frac{1}{\sqrt{2}}   \frac{1}{1 - \frac{1}{2}\lambda^1-\frac{1}{2}\lambda^2} \begin{pmatrix} 1 - \la^2 \\ \la^1 -1\end{pmatrix}.
\]
One can show that
$u_\la$ is bounded as $\la\stackrel{\rm nt}{\to} (1,1)$ 
by a direct calculation.

Alternatively, we can use Proposition \ref{Bptrealzn}. Let
 \[
 e_\pm = \frac{1}{\sqrt{2}}\begin{pmatrix}1 \\ \pm 1 \end{pmatrix}.
 \] 
A realization of the model $(\C^2, u)$ is $(0, e_-, e_-, D)$ where
 $D = e_+ \otimes e_+$, that is,
 \[
D = \frac{1}{2} \begin{bmatrix}1 & 1 \\ 1 & 1 \end{bmatrix}. 
\]
Indeed, on letting $\gamma = e_-$, we have
\begin{align*} 
 (1-D\lambda)^{-1} \gamma  &= \frac{1}{1 - \frac{1}{2}\lambda^1-\frac{1}{2}\lambda^2} \begin{bmatrix} 1 - \frac{1}{2}\lambda^2 & \frac{1}{2}\lambda^2  \\ \frac{1}{2}\lambda^1& 1 - \frac{1}{2}\lambda^1 \end{bmatrix}  \frac{1}{\sqrt{2}} \begin{pmatrix}1 \\ -1 \end{pmatrix} \\ 
&=  \frac{1}{\sqrt{2}}   \frac{1}{1 - \frac{1}{2}\lambda^1-\frac{1}{2}\lambda^2} \begin{pmatrix} 1 - \la^2 \\ \la^1 -1\end{pmatrix} \\
  &= u_\la, \\
  \langle \la u_\la, e_- \rangle &= \frac{1}{\sqrt{2}}   \frac{1}{1 - \frac{1}{2}\lambda^1-\frac{1}{2}\lambda^2} \left\langle \begin{pmatrix} \la^1-\la^1\la^2 \\ \la^1\la^2-\la^2 \end{pmatrix}, \frac{1}{\sqrt{2}} \begin{pmatrix}1 \\ -1 \end{pmatrix}  \right\rangle\\ 
  &= \ph(\la)
 \end{align*}
as required.  Since $(1-D)e_-=e_-=\ga$, we have $\ga\in\mathrm{ran}~(1-D)$, and so $(1,1)$ is a $B$-point for $\ph$ by Proposition \ref{Bptrealzn}.
 \hfill $\Box$

It is a simple matter to modify the foregoing calculations to obtain the following slightly more general example.
\begin{proposition} \label{genex}
If $a_1,a_2\in \C\setminus \{0\}$ then the variety
\[
a_1 + a_2 -a_1\la^1-a_2\la^2 =0
\]
is disjoint from $\D^2$ if and only if $\mathrm{arg}~a_1 = \mathrm{arg}~a_2$.

Let $a_1>0, a_2>0$ and $a_1+a_2=1$.
\begin{enumerate}
\item[\rm (1)] The function
\beq \label{defphi}
\ph(\la) = \frac{a_2\la^1+a_1\la^2 - \la^1\la^2}{1-a_1\la^1-a_2\la^2}
\eeq
is inner on $\D^2$ and has a singularity at $(1,1)$.
\item[\rm(2)]  Let
\[
e_+=\begin{pmatrix} \sqrt{a_1} \\ \sqrt{a_2} \end{pmatrix}, \qquad  e_-= \begin{pmatrix} \sqrt{a_2} \\ -\sqrt{a_1} \end{pmatrix}.
\]

A model for  $\ph$ is $(\C^2,u)$ where
\[
u_\la=\frac{\sqrt{a_1 a_2}(\la^1 - \la^2)}{1-a_1\la^1 -a_2\la^2}e_+ + e_-.
\]
A realization of $(\C^2,u)$ is 
\[
\ph(\la)= \  \langle  \la(1-D\la)^{-1}e_-,e_-\rangle 
\]
where $D=e_+\otimes e_+$.  

\item[\rm(3)]  $(1,1)$ is a $B$-point of $\ph$ and $\ph(1,1) =1$.
\item[\rm(4)]  If $\delta=(\delta^1,\delta^2)$ where $\re \delta^1 > 0, \  \re \delta^2>0$ then the directional derivative $D_{-\delta}\ph(1,1) $ exists and 
\[
D_{-\delta}\ph(1,1) =-\frac{\delta^1\delta^2}{a_1\delta^1+a_2\delta^2}.
\]
\end{enumerate}
\end{proposition}

We leave the proof of the following as an entertaining exercise for the reader.
\begin{proposition}\label{3.2}
For $\phi$ given by equation {\rm (\ref{defphi})} and its model $(\C^2,u)$ as in Proposition {\rm \ref{genex}},  $u_{(1,1)} = e_-$ and 
 \[
Y_{(1,1)} = \left\{ z e_+ + e_-: \mbox{ either } \im z \neq 0 \mbox{ or } z \in \R, \,  -\sqrt{\frac{a_2}{a_1}} < z < \sqrt{\frac{a_1}{a_2}}\right \}.
 \]
 \end{proposition}

The inner function $\ph$ of equation (\ref{eqd8}) illustrates the
complicated nature of the singularities that can occur at $B$-points.  For
any $\al\in\D$ consider the variety $\{\la:\ph(\la)=\al\}$, or, more
precisely,
\[
V_\al \df \{\la\in\C^2: \half \la^1 + \half \la^2 - \la^1\la^2
=\al(1-\half\la^1 - \half\la^2)\}.
\]
It is an exercise to show that $V_\al$ has non-empty intersection with
$\D^2$ and is irreducible.  Clearly $(1,1) \in V_\al$.  It follows that
the cluster set of $\ph$ at $(1,1)$ contains $\D$.  Since $|\ph| \leq 1$
on $\D^2$, the cluster set is $\D^-$.  It is thus impossible to extend
$\ph$ to a continuous function on $\D^2 \cup \{(1,1)\}$.  Nevertheless,
since $(1,1)$ is a $B$-point for $\ph$, the expression $\ph(1,1)$ has a
meaning, and in fact
\[
\ph(1,1) = \lim_{r \to 1-} \ph(r,r) = 1.
\]
It follows that the sets $V_\al \cap \D^2,\, \al\in\D$, all approach $(1,1)$ in
a tangential manner.

It is also interesting to observe how $\ph(\tau)$ behaves as $\tau$
approaches $(1,1)$ along a $C^2$ curve lying in $\T^2$:  $\ph(\tau)$ tends
to a limit $\omega \in \T$, and moreover $\omega$ can be made to equal any
complex number of unit modulus by a suitable choice of $C^2$ curve through
$(1,1)$.

\section{Differentiability at $B$-points} \label{directional}
A function in the Schur class of the bidisk, even if rational, can have singularities at $B$-points on the torus, that is, points to which it cannot be continuously extended, as the foregoing example shows.  It is a remarkable fact that there is nevertheless a rich differentiable structure at such points, as described in Theorem \ref{pickfn}.  The present section is devoted to the proof of this theorem.
First we study limits of $u_\la$ along rays through a $B$-point.

We shall need the notion of the {\em nontangential cluster set} $X_\tau$ of a model $(\M,u)$ of $\ph\in\Schur$ at a $B$-point $\tau\in\partial (\D^2)$.  This is defined to be the set of limits of (weakly or strongly) convergent sequences $\{u_{\la_n}\}$ as $\{\la_n\}$ ranges over all sequences in $\D^2$ that converge nontangentially to $\tau$.  By Corollary \ref{2.10}, for any such sequence  $\{\la_n\}$,  $(1-|\ph(\ln)|)/(1-||\ln||)$ is bounded; it follows that $X_\tau \subset Y_\tau$,
where $Y_\tau$ is defined at the beginning of Section~\ref{JuliaLem}.

\begin{theorem} \label{xdelta}
Let $\tau\in \partial (\D^2)$ be a $B$-point of $\ph\in\Schur$ and let $(\M,u)$ be a model of $\ph$.   For any $\delta\in\Htau$ the nontangential limit (in the norm of $\M$)
\beq\label{defxtau}
x_\tau(\delta) \df  \lim_{\tau-z\de \stackrel {\rm nt}{\to} \tau}u_{\tau-z\delta}
\eeq
 exists in $\M$.  Moreover,
\begin{enumerate}
\item[\rm(1)] $x_\tau(\cdot)$ is a holomorphic $\M$-valued function on $\Htau$;
\item[\rm(2)] $x_\tau(\delta) \in X_\tau$ for all $\delta\in\Htau$;
\item[\rm(3)] $x_\tau(\cdot)$ is homogeneous of degree $0$, that is, $x_\tau(z\delta) = x_\tau(\delta)$ for all $z\in\C$ such that $\delta, z\delta \in \Htau$.
\end{enumerate}
\end{theorem}
Let us spell out the precise meaning of the nontangential limit in equation (\ref{defxtau}).  It is: for any set $S\subset \D^2$ that approaches $\tau$ nontangentially and contains points of the form $\tau - z\de, \ z\in\C$, that tend to $\tau$, the limit of $u_{\tau-z\de}$ as $\tau-z\de \to \tau$ along $S$ exists and equals $x_\tau(\de)$.  In principle the limit depends on $S$, but since the union of two such sets $S$ is a set with the same properties, it is easy to see that the limit is independent of the choice of $S$.
Moreover, for $\delta\in \Htau$, the set $\{\tau-t\de: 0< t <\ep\}$ is contained in $ \D^2$ for suitable $\ep>0$ and approaches $\tau$ nontangentially, so that  there do exist sets $S$ with the required property.
\begin{proof}
We will need the following simple observation.
\begin{lemma}\label{limulan}
Let $(a,\beta,\ga,D)$ be a realization of a model $(\M,u)$ of a function $\ph\in\Schur$.  If $\tau \in\tb$ is a $B$-point of $\ph$ and $x\in Y_\tau$ then, for any sequence $\{\la_n\}$ in $\D^2$,
\[
u_{\la_n} \to x \mbox{  if and only if  }\quad (1-D\la_n)^{-1}D(\la_n-\tau)x \to 0.
\]
\end{lemma}
\begin{proof}
Since $x\in Y_\tau$ equation (\ref{eqaa6}) implies that $(1-D\tau)x =\ga$.  Hence
\begin{eqnarray*}
(1-D\la_n)u_{\la_n} &=& \ga=(1-D\tau)x \\
    &=& (1-D\la_n)x + D(\la_n-\tau)x,
\end{eqnarray*}
and therefore
\[
u_{\la_n} = x + (1-D\la_n)^{-1}D(\la_n-\tau)x,
\]
from which equation the statement follows.
\end{proof}

To resume the proof of Theorem \ref{xdelta}, choose a realization $(a,\beta,\ga,D)$ of the model $(\M,u)$.
We shall prove that the limit in equation (\ref{defxtau}) exists with the aid of a simple identity.  For any $\la,\mu \in \D^2$,
\begin{align}\label{ident}
(1-D\mu)^{-1}D(\mu-\tau)&=[(1-D\mu)^{-1}-(1-D\la)^{-1}]D(\mu-\tau) \nn  \\ &\qquad +(1-D\la)^{-1}D[\mu-\tau-(\la-\tau)]+(1-D\la)^{-1}D(\la-\tau) \nn \\
 &=[(1-D\mu)^{-1}D(\mu-\la)(1-D\la)^{-1}]D(\mu-\tau) \nn \\
 & \qquad + (1-D\la)^{-1} D(\mu-\la) +(1-D\la)^{-1}D(\la-\tau).
\end{align}

  Let $S$ approach $\tau$ nontangentially.  By Proposition \ref{2.9},  $\{u_\la: \la\in S\}$ is bounded in $\M$.  Suppose that $\lim_{z \to 0} u_{\tau-z\de}$  along $\tau-z\de \in S$ does not exist: then there exist sequences $\{w_n\}, \, \{z_n\} \subset \C$ tending to zero and distinct vectors $x,y\in X_\tau$ such that $\tau-z_n\de, \tau-w_n\de \in S$ for all $n$ and
\[
u_{\la_n} \to x,   \qquad u_{\mu_n} \to y
\]
where $\la_n=\tau - z_n\de, \, \mu_n=\tau - w_n\de$.  By 
Lemma \ref{limulan},  
\beq\label {limlan}
(1-D\la_n)^{-1}D(\la_n-\tau)x \to 0.
\eeq
Apply both sides of the identity (\ref{ident}) to the vector $x$ with $\la=\la_n,\, \mu=\mu_n$ to obtain
\begin{align}\label{ABC}
(1-D\mu_n)^{-1}D(\mu_n-\tau )x&=  (1-D\mu_n)^{-1}D(\mu_n-\la_n)(1-D\la_n)^{-1}D(\mu_n-\tau)x  \nn \\
 & \qquad + (1-D\la_n)^{-1} D(\mu_n-\la_n)x \nn \\ &\qquad \qquad+(1-D\la_n)^{-1}D(\la_n-\tau)x\nn  \\
  &= A_n + B_n + C_n, \mbox{  say}.
\end{align}
We shall show that the right hand side of this equation tends to zero.  By the relation (\ref{limlan}), $C_n\to 0$.  Furthermore, since
\[
 \mu_n-\la_n = (z_n-w_n)\de =\frac{w_n-z_n}{w_n}(\mu_n-\tau),
\]
and
\[
\mu_n-\tau = -w_n\de =\frac{w_n}{z_n}(\la_n-\tau),
\]
we have
\begin{align*}
A_n &=(1-D\mu_n)^{-1}D\frac{w_n-z_n}{w_n}(\mu_n-\tau)(1-D\la_n)^{-1}D\frac{w_n}{z_n}(\la_n-\tau)x\\
   &= \frac{w_n-z_n}{z_n}(1-D\mu_n)^{-1}D(\mu_n-\tau)(1-D\la_n)^{-1}D(\la_n-\tau)x.
\end{align*}
Here the operators $(1-D\mu_n)^{-1}D(\mu_n-\tau)$ are uniformly bounded, since $\mu_n \to \tau$ nontangentially, while the vectors $(1-D\la_n)^{-1}D(\la_n-\tau)x$ tend to zero by (\ref{limlan}).  
By passing to a subsequence of $\{\mu_n\}$ if necessary, we may ensure that $w_n/z_n$ is bounded.  Thus $A_n \to 0$ as $n \to \infty$.

Similarly, since
\[
 \mu_n-\la_n = (z_n-w_n)\de =\frac{w_n-z_n}{z_n}(\la_n-\tau),
\]
we have
\[
B_n=(1-D\la_n)^{-1}D(\mu_n-\la_n)x = \frac{w_n-z_n}{z_n} (1-D\la_n)^{-1}D(\la_n-\tau)x,
\]
and so $B_n \to 0$ as $n\to \infty$.
It follows from equation (\ref{ABC}) that
\[
(1-D\mu_n)^{-1}D(\mu_n-\tau)x \to  0,
\]
and therefore, by Lemma \ref{limulan}, that $u_{\mu_n} \to x$, a contradiction.  Thus the limit 
\beq\label{defxS}
x_\tau(\delta)=\lim_{z\to 0} u_{\tau-z\de} \mbox{  along  } \tau-z\de \in S
\eeq
exists in $\M$ as required.\\

\noindent (1)  Define
\[
F_t(\delta)= u_{\tau-t\de}
\]
whenever $\tau-t\de\in\D^2$.  Observe that $\tau-t\de\in\D^2$ whenever $0 < t < r(\delta)$, where, for $\de\in\Htau$,
\[
r(\de) \df  \left\{ 
\begin{array}{lcl} \ds 2\min_{j=1,2}  \frac{2\re \overline{\tau^j}\de^j}{|\de^j|^2} & \mbox{ if } & \tau\in\T^2, \\
 \ds  \min\left\{ \frac{2\re \overline{\tau^1}\de^1}{|\de^1|^2}, \frac{1-|\tau^2|}{|\de^2|} \right\} &  \mbox{ if } & \tau\in\T\times \D. \end{array} \right.
\]
Consider any $\delta_0\in\Htau$ and pick a number $c>0$ and a compact neighborhood $U$ of $\delta_0$ such that $r(\delta)>\half r(\delta_0)$ and $\min_{|\tau_j|=1}  \re \overline{\tau^j}\delta^j > c$ for all $\delta\in U$.  Then $F_t$ is defined and holomorphic on the interior of $U$ for all $t \in (0,\half r(\delta_0))$.  

We claim that the set
\[
 S=\{\tau -t\de: \delta\in U, \, 0 < t < \tfrac 14 r(\delta_0) \}
\]
approaches $\tau$ nontangentially.  Consider $\tau\in\T^2$.   If $\delta\in U$ and $0< t< \tfrac 14 r(\delta_0)$ then
\begin{align*}
\frac{||t\de||}{1-||\tau-t\de||} &= \frac{t||\delta||(1+||1-t\tau^*\delta||)}{1-||1-t\tau^*\delta ||^2} \\
  & \leq \frac{2t||\delta||}{1-\max_j(1-2t\re \overline{\tau^j}\delta^j + t^2|\delta^j|^2)} = \frac {2||\delta||}{\min_j (2\re\overline{\tau^j}\delta^j - t|\delta^j|^2)}.
\end{align*}
Now, since $|\delta^j| > c$,
\begin{align*}
\min_j (2\re\overline{\tau^j}\delta^j - t|\delta^j|^2) &= \min_j |\delta^j|^2\left(\frac{2\re\overline{\tau^j}\delta^j}{|\delta^j|^2} - t\right)\\
   & > c^2 \min_j \left(\frac{2\re\overline{\tau^j}\delta^j}{|\delta^j|^2} -t\right) = c^2(r(\delta)-t)\\
  & > c^2 (\half r(\delta_0) - \tfrac 14 r(\delta_0)) = c^2 \tfrac 14 r(\delta_0).
\end{align*}
Hence 
\[
\frac{||\tau - (\tau -t\de)||}{1-||\tau-t\de||} < \frac{8||\delta||}{c^2r(\delta_0)},
\]
which is bounded.    Alternatively, suppose $\tau\in\T\times \D$.   For small enough $t>0$ we have $||\tau -t\de|| = |\tau^1-t\de^1|$ and so
\begin{align*}
1-\|\tau-t\de\|^2 &=  2t\re (\overline{\tau^1}\de^1)- t^2|\de^1|^2 \\
  &>  t c^2\tfrac 14 r(\de_0).
\end{align*}
Hence, for small $t>0$,
\begin{align*}
\frac{||\tau - (\tau -t\de)||}{1-||\tau-t\de||} &= \frac{||\de||(1+||\tau-t\de||)}{2\re \de^1- t|\de^1|^2} \\
   &<\frac{4C\|\de\|}{c^2r(\de_0)},
\end{align*}
where $C$ is an upper bound for $1+||\tau-t\de||$ over $t\in[0, \tfrac 14], \ \de\in U$. Thus, in either case, $S$ approaches $\tau$ nontangentially.

Since $\tau$ is a $B$-point of $\ph$ and $S$ approaches $\tau$ nontangentially, it follows from Proposition \ref{2.9} that $\{u_\la:\la\in S\}$ is bounded.  In other words, the holomorphic functions $F_t, \, 0<t< \tfrac 14 r(\delta_0),$ are uniformly bounded on $U$.
Hence their pointwise limit as $t\to 0$, which is $x_\tau(\cdot)$, is holomorphic on $U$.  As $U$ is a neighborhood of a general point of $\Htau$, $x_\tau$ is holomorphic on $\Htau$.  This establishes (1).\\

\noindent (2) is immediate from the definition of the nontangential cluster set $X_\tau$.\\

\noindent (3)  Consider $\de\in\Htau$ and $z\in\C$ such that $z\de\in\Htau$.  For a suitable $\ep>0$ the set 
\[
S=\{\tau-t\de, \tau - tz\de : 0 < t < \ep\}
\]
is contained in $\D^2$ and approaches $\tau$ nontangentially.  For any sequence $\{t_n\}$ tending to $0$ in $(0,\ep)$,  the sequences $\{\tau-t_n\de\}$ and $\{\tau-t_nz\de\}$ are both sequences in $S$ that tend to $\tau$, and hence the limit of $u$ along either of these sequences is $x_\tau(\de)$ by equation (\ref{defxtau}).  Its limit along the second is also $x_\tau(z\de)$.
Thus $x_\tau(\cdot)$ is homogeneous of degree $0$.
\end{proof}
The vectors $x_\tau(\delta)$, defined in (\ref{defxtau})
by
$ \displaystyle
x_\tau(\delta) =  \lim_{\tau-z\de \stackrel {\rm nt}{\to} \tau}u_{\tau-z\delta}
$,
enable one
to calculate the directional derivatives at $\tau$ in the following way.
\begin{theorem}\label{directderiv}
Let $\tau\in\partial(\D^2)$ be a $B$-point of $\ph\in\Schur$.  For any $\de\in\Htau$  the directional derivative
$D_{-\de}\ph(\tau)$ exists and is analytic and homogeneous of degree $1$ in $\de$.  Moreover, for any model $(\M,u)$ of $\ph$ and any $y\in Y_\tau$,
\beq \label{formDph}
D_{-\de}\ph(\tau)= - \ph(\tau) \langle \de x_\tau(\de), \tau y \rangle.
\eeq
In particular,
\beq
\label{eqd7}
D_{-\de}\ph(\tau)= - \ph(\tau) \sum_{|\tau^j|=1} \overline{\tau^j}\de^j||x_\tau^j(\de)||^2.
\eeq
\end{theorem}
\begin{proof}
Let $\ph(\tau)=\omega$ and choose a model $(\M,u)$ of $\ph$.   By Proposition  \ref{getx} we have, for $\la\in\D^2$ and $y\in X_\tau$,
\[
1-\overline{\omega}\phi(\lambda) =  \sum_{|\tau^j|=1}(1 - \overline{\tau^j}\lambda^j) \langle u_\lambda^j, y^j \rangle.
\] 
On substituting $\la=\tau-t\de$ 
and multiplying through by $-\omega$ we have
\[
\ph(\tau-t\de) - \ph(\tau) = -\omega t \langle\de u_{\tau-t\de}, \tau y \rangle.
\]
Divide by $t$ and let $t\to 0+$ to obtain equation (\ref{formDph}).
Since $x_\tau(.)$ is analytic on $\Htau$, so is $D_{-\de}\ph(\tau)$ as a function of $\de$; since $x_\tau(.)$ is homogeneous of degree $0$ in $\de$,  $D_{-\de}\ph(\tau)$ is homogeneous of degree $1$.
Furthermore, since equation (\ref{formDph}) holds for any $y\in Y_\tau$, it holds when $y = x_\tau(\de)$.  Since $y^j=0$ when $|\tau^j|< 1$ we deduce equation (\ref{eqd7}).
\end{proof}

The following observation can be regarded as an analog of equation (\ref{valeta2}) in the classical Julia-Carath\'eodory theorem -- indeed, it can easily be proved by the application of the classical theorem, Theorem \ref{thma1}, to the one-variable function $\psi(\la)\df\ph(\la\tau)$.  We give an independent proof.
\begin{corollary} \label{radialderiv}
If $\tau\in\tb$ is a $B$-point for $\ph\in\Schur$ then
\[
D_{-\tau}\ph(\tau) = -\ph(\tau)\liminf_{\la\to\tau} \frac{1-|\ph(\la)|}{1-\|\la\|}.
\]
\end{corollary}
\begin{proof}
By definition of $x_\tau$ and Theorem \ref{2.16},
\[
x_\tau(\tau)= \lim_{1-z \nt 1} u_{\tau-z\tau} = u_\tau.
\]
Combine this fact with equation (\ref{eqd7}) to obtain
\[
D_{-\tau}\ph(\tau) = -\ph(\tau) \sum_{|\tau^j|=1} \|u_\tau^j\|^2 = -\ph(\tau)\|u_\tau\|^2,
\]
and with equation (\ref{normutau}) to establish the corollary.
\end{proof}

\begin{proposition} \label{posxtau}
Let $\ph\in\Schur$, let $\tau\in\tb$ be a $B$-point of $\ph$ and let $\ph(\tau)=\omega$.  For any model $(\M,u)$, any $y\in Y_\tau$ and any $ \de\in\Htau$,
\beq \label{condxtau}
\re  \langle \de x_\tau(\de), \tau y \rangle \geq 0.
\eeq
\end{proposition}

\begin{proof}  For sufficiently small $t>0$ we have $\tau-t\de \in\D^2$ and so
$|\ph(\tau-t\de)|^2 \leq 1$.  
As
\[
\ph(\tau-t\de) \ = \ 
\omega+tD_{-\de}\ph(\tau) + o(t),
\]
we get that, for small $t>0$,
\[
2t \re (\bar\omega D_{-\de}\ph(\tau)) + o(t) \leq 0,
\]
and so $\re (\bar\omega D_{-\de}\ph(\tau))  \leq 0.$
By equation (\ref{formDph}) we have
\[
\re(\bar\omega(-\omega \langle \de x_\tau(\de), \tau y\rangle)) \leq 0,
\]
which implies the inequality (\ref{condxtau}).
\end{proof}
We come to one of the main results of the paper: the existence of a function in the Pick class associated with any $B$-point of $\ph$ (see Theorem \ref{pickfn}).
\begin{theorem} \label{pickfunc}
Let $\tau\in\T^2$ be a $B$-point for $\ph\in\Schur$.  There exists a function $h$ in the Pick class, analytic and real-valued on $(0,\infty)$ and
with $-zh(z)$ also in the Pick class, and such that 
\beq\label{valh1}
h(1) = -\liminf_{\la\to\tau}\frac{1-|\ph(\la)|}{1-\|\la\|}
\eeq 
and, for all $\de\in\Htau$,
\beq\label{defproph}
D_{-\de}\ph(\tau) = \ph(\tau)\overline{\tau^2}\de^2 h\left(\frac{\overline{\tau^2}\de^2}{\overline{\tau^1}\de^1}\right).
\eeq
\end{theorem}
\begin{proof}
Let $\omega=\ph(\tau)$.  Choose a model of $\ph$ and any $y\in Y_\tau$.
Since $ \de x_\tau(\de)$ is a homogeneous function of degree $1$ in $\de$ on $\Htau$, the function $(\de^1)^{-1}\lan \de x_\tau(\de), \tau y\ran $ is a homogeneous function of degree $0$, hence is a function of $\de^2/\de^1$.  We may therefore define a function $g$ by
\beq \label{defg}
g( \zeta^2/\zeta^1) = (\zeta^1)^{-1}\lan \de x_\tau(\de), \tau y\ran
\eeq
where $\zeta\in\HH$ and $\de=\tau\zeta\in \Htau$.  This function $g$ is analytic on the set
\[
\mathcal{D}(g) \df \left\{\zeta^2/\zeta^1: \zeta \in \HH  \right\} =  \C\setminus (-\infty,0],
\]
and in view of Proposition \ref{posxtau} we have
\begin{align*}
0 &\leq  \re \lan \de x_\tau(\de), \tau y\ran = \re \left( \zeta^1 g( \zeta^2/\zeta^1)\right)
\end{align*}
for all $\zeta\in\HH$.  Thus $g$ satisfies (1) in the following assertion.
\begin{lemma}
Let $g$ be analytic on $\C\setminus (-\infty,0]$.  The following statements are equivalent.
\begin{enumerate}
\item[\rm (1)] $\re (\la g(z)) \geq 0$ for all $z\in \C\setminus (-\infty,0]$ and all $\la\in\Ha$ such that $\la z\in \Ha$;
\item[\rm (2)]
 $g$ is argument-decreasing, in the sense that $0 \leq \arg g(z) \leq \arg z $ whenever $ 0 \leq \arg z < \pi$;
\item[\rm(3)] the functions $g$ and 
\[
h(z) \df  - g(z)/z
\]
are in the Pick class and real on $(0,\infty)$.
\end{enumerate}
\end{lemma}
We leave the proof of this lemma to the reader.
Applying it to the function $g$ in definition (\ref{defg}) we obtain a function $h$ in the Pick class, analytic and real-valued on $(0,\infty)$, and satisfying $g(z)=-zh(z)$ for $z\in\Pi$.  By equations (\ref{formDph}) and (\ref{defg}) we have, if $\zeta=\tau^*\de$,
\begin{align*}
D_{-\de}\ph(\tau) &= -\omega \lan \de x_\tau(\de),\tau y\ran \\
  &= -\omega\zeta^1g(\zeta^2/\zeta^1) \\
  &= \omega\zeta^2 h(\zeta^2/\zeta^1)\\
  &= \omega\overline{\tau^2}\de^2 h\left(\frac{\overline{\tau^2}\de^2}{\overline{\tau^1}\de^1}\right)
\end{align*}
as required.

On combining Corollary \ref{radialderiv} with equation (\ref{defproph}) we have
\[
  -\ph(\tau)\liminf_{\la\to\tau} \frac{1-|\ph(\la)|}{1-\|\la\|} =D_{-\tau}\ph(\tau) = \ph(\tau) h(1)
\]
and so equation (\ref{valh1}) is true.
\end{proof}
Here is the corresponding statement for $B$-points in $\tb\setminus \T^2$.
\begin{proposition} \label{Dphundist}
If $\tau\in\T\times\D$ is a $B$-point for $\ph\in\Schur$ then, for any $\de\in\tau^1\Ha\times\C$,
\[
D_{-\de}\ph(\tau) = -\ph(\tau)\overline{\tau^1}\de^1 \liminf_{\la\to\tau} \frac{1-|\ph(\la)|}{1-\|\la\|}.
\]
\end{proposition}
This formula corresponds to the choice $h(z)=-\mathrm{const}/z$ in Theorem \ref{pickfunc}.
\begin{proof}
Let $\ph(\tau)=\omega$  and let
\[
\al = \liminf_{\la\to\tau} \frac{1-|\ph(\la)|}{1-\|\la\|}.
\]
  Choose a model of $\ph$ and any $y\in Y_\tau$.  Then $y^2=0$ and, by Theorem \ref{directderiv},
\beq \label{Dde}
D_{-\de}\ph(\tau) = -\omega\de^1\lan x_\tau^1(\de), \tau^1 y^1\ran.
\eeq
Since $x^1$ is analytic on $\tau^1\Ha\times\C$ we may define an entire function $F$ by
\[
F(\zeta) = \lan x_\tau(\tau^1,\zeta), \tau^1 y^1 \ran, \qquad \zeta\in\C.
\]
By the $0$-homogeneity of $x_\tau$ on $\tau^1\Ha\times\C$,
\[
x_\tau^1(z\tau^1 , z\zeta) = x_\tau^1(\tau^1,\zeta)
\]
for all $z\in\Ha$.  Hence, if $\de\in\tau^1\Ha\times\C$ then $z\df \overline{\tau^1}\de^1$ is in $\Ha$, and
\begin{align} \label{DdeF}
\lan \de^1 x_\tau^1(\de), \tau^1 y^1 \ran &= \de^1\lan x_\tau^1(z\tau^1, z\tau^1\de^2/\de^1), \tau^1 y^1 \ran \nn \\
  &= \de^1 \lan x_\tau^1(\tau^1, \tau^1\de^2/\de^1), \tau^1 y^1 \ran \nn \\
 &= \de^1 F(\tau^1\de^2/\de^1).
\end{align}
By Proposition \ref{posxtau}, for all $\de\in \tau^1\Ha\times\C$,
\beq \label{strongineq}
\re \left( \de^1F(\tau^1\de^2/\de^1)\right) \geq 0
\eeq
and in particular, when $\de^1=\tau^1$, for all $\de^2\in\C$,
\[
\re \left(\tau^1 F(\de^2)\right) \geq 0.
\]
It follows from Liouville's Theorem that $F$ is constant -- say $F=\overline{\tau^1}c$,  so that equations (\ref{Dde}) and (\ref{DdeF}) yield, for all $\de\in\tau^1\Ha\times\Ha$,
\[
D_{-\de}\ph(\tau) = -\ph(\tau) \de^1\overline{\tau^1} c.
\]
In particular, when $\de=\tau$, this relation and Corollary \ref{radialderiv} combine to give
\[
-\ph(\tau) \al = D_{-\tau}\ph(\tau) = -\ph(\tau) c,
\]
so that $c=\al$ and the proposition is proved.
\end{proof}

\section{ Models and $C$-points}\label{holdiffC}

In this section we show that the function-theoretic property of being a $C$-point of $\ph\in\Schur$ corresponds to a natural property of models of $\ph$: roughly, that $u_\la$ extends continuously to the point in question on nontangential approach sets.
Recall Definition \ref{4.4}: $\tau\in\partial (\D^2)$ is a $C$-point for $\ph\in\Schur$ if 
$\ph$ has a holomorphic differential on every set that approaches $\tau$ nontangentially, and if, furthermore, $|\ph(\tau)|=1$.

We first show that every $C$-point is a $B$-point.
\begin{proposition} \label{4.7}
Let $\phi \in \mathcal{S}, \tau \in \partial (\D^2)$ and let $S \subset \mathbb{D}^2$ approach $\tau$ nontangentially.  
If
\[
\lim_{\la\to\tau,\, \, \la\in S} \ph(\la) = \omega \in\T
\]
and $\ph$ has a holomorphic differential at $\tau$ on $S$ then $\tau$ is a $B$-point for $\ph$.
\end{proposition}

\begin{proof} 
By hypothesis there exist $ \eta^1, \eta^2 \in \mathbb{C}$ such that the equations {\rm (\ref{4.5})} and {\rm (\ref{4.6})} hold. 
For $\lambda \in S$, 
\begin{align} \label{4.8}
1 - |\phi(\lambda)| & \le 1 - |\phi(\lambda)|^2  = 1-|\omega+\eta\cdot(\la-\tau)+o(\|\tau-\la\|)|^2  \nn \\
    & = 2 \mathrm{Re}(\overline{\omega}\eta^1(\tau^1-\lambda^1)+\overline{\omega}\eta^2(\tau^2-\lambda^2))+o(||\tau-\lambda||) \nn \\
 &\leq 2 ||\eta|| \cdot ||\tau-\lambda|| + o(||\tau-\lambda||).
 \end{align}
Since $S$ approaches $\tau$ nontangentially there is a sequence $\{\lambda_n\} \subset S$ and a constant $c$ such that $\lambda_n \rightarrow \tau$ and
\[ 
||\tau - \lambda_n|| \le c (1-||\la_n||).
 \]
Hence  (\ref{4.8}) implies that 
$\frac{1 - |\phi(\lambda_n)|}{1-||\la_n||}$  is  bounded.
\end{proof}

Consider $\tau \in \partial (\D^2)$ and let $S \subset \mathbb{D}^2$  approach $\tau$ nontangentially. If $\delta \in \mathbb{C}^2$, let us say that $S$ {\em admits} $\delta$ if there exists $\epsilon >0$ such that
 \[
 \tau -t\delta \in S \mbox{ for } t \in (0, \epsilon).
 \]
More generally, if $\Delta \subset \mathbb{C}^2$ we say that $S$ admits $\Delta$ if there exists $\epsilon >0$ such that 
\[ 
\tau-t\delta \in S \mathrm{ \ whenever \ } t \in (0, \epsilon) \mathrm{ \ and \ } \delta \in \Delta. 
\]
We say that $S$ is {\em untapered at} $\tau$ if $S$ approaches $\tau$ nontangentially and admits a nonempty open set.

\begin{theorem}\label{4.19}
Let  $\phi\in\Schur$, $\tau\in \partial (\D^2)$ and let  $S \subset \mathbb{D}^2$ be untapered at  $\tau$.
 The following three conditions are equivalent. 
\begin{enumerate}
\item[\rm(1)]
\blue
$\lim_{\la \to\tau,\la\in S} |\phi(\la)|=1$ and 
\black
 for every model  $(\mathcal{M}, u)$ of $ \phi$,  $u_\lambda$ extends  continuously  to $S^-$; 

\item[\rm(2)] \blue $\lim_{\la \to\tau,\la\in S} |\phi(\la)|=1$ and \black there exists a model  $(\mathcal{M}, u)$ of $ \phi$ such that  $u_\lambda$ extends  continuously  to $S^-$; 

\item[\rm(3)]  $ \phi $ has  a  holomorphic  differential  at $ \tau$ on $ S$  with  $|\ph(\tau)| = 1$.
\end{enumerate} 
\end{theorem}
\begin{proof}  Clearly (1) implies (2). We shall show that (2) implies (3). Accordingly, fix a model $(\mathcal{M},u)$ of $\phi$  such that $u_\lambda$ extends  continuously to $S^-$ -- say $u_\lambda \rightarrow x$ as $\lambda \rightarrow \tau$ in $S$. 
\blue
Hence $u_\la$ is bounded on some sequence that tends to $\tau$ nontangentially.
By Corollary \ref{2.10}, $\tau$ is a B-point of $\ph$.

\black

By Proposition \ref{getx},  $x^j=0$ if $|\tau^j|<1$ and there exists $\omega\in\T$ such that, for all $\lambda \in \D^2$,
\begin{align*} 
1 - \overline{\omega}\phi(\lambda) & = \sum_{|\tau^j|=1} (1 - \overline{\tau^j}\lambda^j)\langle u_\lambda^j, x^j \rangle  \\ 
 &= \sum_{|\tau^j|=1} (1 - \overline{\tau^j}\lambda^j) || x^j ||^2 + \sum_{|\tau^j|=1} (1 - \overline{\tau^j}\lambda^j)\langle u_\lambda^j - x^j, x^j \rangle.
\end{align*}
Hence, for all $\la\in S$,
\[
\phi(\lambda) = \omega+\omega \sum_{|\tau^j|=1}\overline{\tau^j} ||x^j||^2 (\lambda^j-\tau^j)  + e(\lambda) 
\]
where
\[ 
e(\lambda) = - \omega \sum_{|\tau^j|=1} (1 - \overline{\tau^j}\lambda^j)\langle u_\lambda^j -x^j, x^j\rangle.
 \]
 By the Cauchy-Schwarz inequality,
\[ 
|e(\lambda)| \le ||\lambda - \tau||\cdot ||u_\lambda -x||\cdot ||x||. 
\]
Since $u_\lambda \rightarrow x$ as $\lambda \rightarrow \tau$ in $S$, this shows that $e(\la)=o(||\la-\tau||)$.
Hence
\[
\ph(\la)=\omega + \eta^1(\la^1-\tau^1) + \eta^2(\la^2-\tau^2) + o(||\la-\tau||)
\]
for $\la\in S$, where $\eta^j= \omega \overline {\tau^j} ||x^j||^2$ for $ j=1,2$.
Thus $(3)$ holds. \\

\noindent (3)$\Rightarrow$(1)  will follow from three lemmas.
 \begin{lemma} \label{constxdel}
If $\ph\in\Schur$ has a holomorphic differential on a set $S$ untapered at $\tau\in\tb$ and $|\ph(\tau)|=1$ then, for any model $(\M,u)$ of $\ph$,  the function $x_\tau(\cdot)$ is constant on $\Htau$ and equal to $u_\tau$, where $u_\tau$ is defined in Theorem {\rm \ref{2.16}}.  Furthermore the holomorphic differential of $\ph$ on $S$ is
\beq \label{diffonS}
\ph(\la) = \ph(\tau)\left(1 + \overline{\tau_1}||u_\tau^1||^2 (\la^1 -\tau^1)+  \overline{\tau_2}||u_\tau^2||^2 (\la^2-\tau^2)\right) + o(||\la-\tau||).
\eeq
\end{lemma}
\begin{proof}

Let $(\M,u)$ be any model of $\ph$.  Since $S$ is untapered at $\tau$ there is an open set $\Delta$ in $\C^2$ such that $S$ admits $\Delta$.  By hypothesis there exist $\omega\in \T$ and $\eta^1,\eta^2\in\R$ such that, for $\la\in S$, 
\[
\ph(\la)=\omega +\eta^1(\la^1-\tau^1)+\eta^2(\la^2-\tau^2)+e(\la),
\]
where
\[
\lim_{\la\to\tau, \la\in S} \frac{e(\la)}{\|\la-\tau\|}=0.
\]
In particular, for $\de\in\Delta$ and $t>0$ such that $\tau-t\de \in S$,
\begin{equation*}\tag{8.5a} \label{8.5a}
\ph(\tau-t\de) = \omega -t \eta^1\de^1 -t\eta^2\de^2 + o(t).
\end{equation*}
On the other hand, by Theorems \blue \ref{xdelta} and \black\ref{directderiv}  we have, for $\de\in\Htau$,
\[
D_{-\de}\ph(\tau) = -\omega\langle \de x_\tau(\de), \tau x_\tau(\de) \rangle,
\]
\blue where $x_\tau(\cdot)$ is a holomorphic $\M$-valued function on $\Htau$,\black~ 
and so
\begin{equation*}\tag{8.5b} \label{8.5b}
\ph(\tau-t\de) = \ph(\tau)+tD_{-\de}\ph(\tau) +o(t)=\omega -t\omega \langle \de x_\tau(\de),\tau  x_\tau(\de) \rangle +o(t).
\end{equation*}

Comparison of the linear terms in these two approximations for $\ph(\tau-t\de) $ shows that, for $\de\in \Delta$,
\begin{equation*}\tag{8.5c}\label{8.5c}
\eta^1\de^1 +\eta^2\de^2 = \omega \langle \de x_\tau(\de), \tau x_\tau(\de) \rangle = \omega \overline{\tau^1}\de^1 ||x_\tau^1(\de)||^2 + \omega \overline{\tau^2}\de^2||x_\tau^2(\de)||^2.
\end{equation*}
 \blue
Rearranging, we get
\begin{equation*}\tag{8.5d}\label{eq84}
\delta^1( \eta^1 - \tau^1 \| x^1(\delta)\|^2) + \delta^2( \eta^2 - \tau^2 \| x^2(\delta)\|^2) = 0.
\end{equation*}
Initially, we know \eqref{eq84} holds on $\Delta$, but by Theorem \ref{xdelta}  the left-hand side is real analytic in $\delta$   for $\delta$ in  ${\mathbb H}(\tau)$, so since a non-zero real analytic function cannot vanish on a non-empty open subset of a connected set, we know that \eqref{eq84} holds on all of ${\mathbb H}(\tau)$, and so there we have
\begin{equation*}\tag{8.5e}\label{8.5e}
\delta^1 / \delta^2 = - ( \eta^2 - \tau^2 \| x^2(\delta)\|^2)/ ( \eta^1 - \tau^1 \| x^1(\delta)\|^2).
\end{equation*}
If we knew $\eta^1$ and $\eta^2$ were real, we could immediately conclude that 
$ \eta^2 - \tau^2 \| x^2(\delta)\|^2=0 $ and $ \eta^1 - \tau^1 \| x^1(\delta)\|^2 =0$
since the ratio $\delta^1/\delta^2$ contains an open set, and therefore cannot always be real.

Otherwise, take $\delta \in {\mathbb H}(\tau)$ real, which we can do since we have assumed $\tau$ is real, and take the imaginary part of \eqref{eq84}.  Then we get
$$\delta^1 \Im (\eta^1) + \delta^2 \Im (\eta^2) = 0.$$
Perturbing $\delta^1$ and $\delta^2$ while keeping them real gives $\Im(\eta^1) = 0 = \Im(\eta^2)$.
Since $\eta^1$ and $\eta^2$ are now real, we can apply the previous argument.
Therefore we deduce that
\black
\beq \label{valeta}
\eta^1 = \overline{\tau^1}\omega||x^1_\tau(\de)||^2, \qquad  \eta^2=\omega\overline{\tau^2} ||x^2_\tau(\de)||^2,
\eeq
and so 
\beq\label{constnorm}
||x_\tau(\de)||^2 = \bar\omega\eta^1\tau^1+\bar\omega\eta^2\tau^2
\eeq
for all $\de\in\Delta$ (recall that, by Proposition \ref{getx},  $\eta^j=||x^j_\tau(\de)||^2=0$ if $|\tau^j|<1$).

By Theorem \ref{xdelta}, $x_\tau$ is an analytic $\M$-valued function on $\Htau$, and we have shown that $||x_\tau(\cdot)||$ is constant on $\Delta$.   It follows that $x_\tau$  is constant on each connected component of $\Delta$.  Thus $x_\tau$ is a holomorphic function on the connected set $\Htau$ and is constant on a nonempty open set.  Hence $x_\tau$ is constant on $\Htau$.   

As $t\to 0$ in the interval $(0,1]$ we have $\tau-t\tau \stackrel{\rm nt}{\to}\tau$, and so, by equation (\ref{defxtau}),
\[
x_\tau(\tau) = \lim_{t\to 0+} u_{\tau-t\tau}.
\]
Theorem \ref{2.16} implies that the last limit is $ u_\tau$, and hence the constant value of $x_\tau$ is $u_\tau$.   It follows from equations (\ref{valeta}) that
\[
\eta^1 = \overline{\tau^1}\omega||u_\tau^1||^2, \qquad  \eta^2=\omega\overline{\tau^2} ||u^2_\tau||^2,
\]
 and so equation (\ref{diffonS}) holds on $S$.  
\end{proof}
\begin{lemma}\label{xconstCpt}
If $(\M,u)$ is a model of $\ph\in\Schur$ and $x_\tau(\cdot)$ is a constant function for some $B$-point $\tau\in\T^2$ of $\ph$ then the value of $x_\tau$ is $u_\tau$. Moreover, for any  set $S$ that approaches $\tau$ nontangentially, $u_\la$ extends continuously to $S^-$ and
\[
\lim_{\la\stackrel{\rm nt}{\to} \tau} u_\la = u_\tau.
\]

\end{lemma}
\begin{proof}
As in the preceding Lemma, the constant value of $x_\tau$ is $x_\tau(\tau) = u_\tau$, by Theorem \ref{2.16}.    We wish to show that $u_\la \to u_\tau$ as $\la \to \tau$ nontangentially.   Consider any sequence $\{\la_n\}$  converging nt to $\tau$.   Let
\[
t_n = ||\tau - \la_n||, \quad \de_n= \frac{\tau-\la_n}{||\tau-\la_n||},
\]
so that $t_n>0, \, ||\de_n|| = 1$ and $\la_n= \tau-t_n\de_n \in \D^2$.  
If the aperture of $\{\la_n\}$ is $c>0$, then
\begin{align*}
c & \geq \frac{||\tau-\la_n||}{1-||\la_n||} = \frac{t_n(1+||\la_n||)}{1-||\tau-t_n\de_n||^2}\\
  & > \frac{t_n}{t_n\min_j (\re 2\overline{\tau^j}\de_n^j - t_n|\de_n^j|^2)} > \frac{1}{2\min_j \re \overline{\tau^j}\de_n^j}.
\end{align*}
Thus $\re \overline{\tau^j}\de_n^j > (2c)^{-1} > 0$ for all $n$ and $j=1,2$.

Let $\de_\infty$ be a cluster point of $\{\de_n:n\geq 1\}$.  By passing to a subsequence of $\{\la_n\}$ we can suppose that $\de_n \to \de_\infty$.  We have $||\de_\infty||=1$ and $\re \overline{\tau^j}\de_\infty^j \geq (2c)^{-1}$ for $j=1,2$, so that $\de_\infty\in\Htau$.     Let $\mu_n=\tau-t_n \de_\infty$. Then by the definition (\ref{defxtau}) of $x_\tau$ and the constancy of $x_\tau$ we have 
\[
u_{\mu_n} = u_{\tau-t_n\de_\infty} \to x_\tau (\de_\infty) = u_\tau \mbox{  as  } n\to \infty.
\]
Choose any realization $(a,\beta,\ga,D)$ of $(\M,u)$.  Then $u_\la= (1-D\la)^{-1}\ga$ for every $\la\in\D^2$, and so
\begin{align} \label{resolvid}
u_{\la_n} - u_{\mu_n} & = [(1-D\la_n)^{-1}  - (1-D\mu_n)^{-1}]\ga \nn\\
   &= (1-D\la_n)^{-1}D(\la_n-\mu_n)(1-D\mu_n)^{-1}\ga \nn\\
  &=(1-D\la_n)^{-1}Dt_n(\de_\infty - \de_n) u_{\mu_n}.
\end{align}

The right hand side converges to zero in $\M$ as $n\to\infty$.  Indeed, the vectors $u_{\mu_n} \to u_\tau$, the operators $\de_\infty - \de_n$ tend to zero in the operator norm,
and the operators $(1-D\la_n)^{-1}Dt_n$ are uniformly bounded: 
\[
||(1-D\la_n)^{-1}D t_n|| 
\ \leq \
 \frac{||\tau -\la_n||}{1-||\la_n||}\ \leq \ c .
\]

We have shown that every sequence $\{\la_n\}$ converging nt to $\tau$ has a subsequence for which the corresponding $u_{\la_n}$ tend to $u_\tau$.  It follows that $u_{\la_n} {\to} u_\tau$ for every sequence $\{\la_n\}$ converging nt to $\tau$.  Hence $u_\la$ extends continuously to $S^-$ for any set $S$ that approaches $\tau$ nontangentially.  
\end{proof}

\begin{lemma} \label{CptTD}
If $(\M,u)$ is a model of $\ph\in\Schur$ and $\tau\in\T\times\D$ is a $B$-point for $\ph$ then, for any set $S$ that approaches $\tau$ nontangentially, $u_\la$ extends continuously to $S^-$ and
\[
\lim_{\la\stackrel{\rm nt}{\to} \tau} u_\la = u_\tau.
\]
\end{lemma}
\begin{proof}
  We claim that $Y_\tau$ (and hence also its nonempty subset $X_\tau$) consists of a single point.  Let $x\in Y_\tau$.  If $\mathcal{N}$ is the closed linear span of  $\{u_\la^1:\la\in\D^2\}$ in $\M^1$ then $x^1\in \mathcal{N}$.  By Proposition \ref{getx} $x^2=0$ and there exists $\omega\in\T$ such that, for all $\la\in\D^2$,
\[
1-\overline\omega\ph(\la) = (1-\overline{\tau^1}\la^1) \lan u_\la^1, x^1 \ran.
\]
Thus $\lan u_\la^1, x^1 \ran$ is uniquely determined for every $\la\in\D^2$.  Since $x\in\mathcal{N}$, it follows that $x$ is uniquely determined.  Thus $X_\tau$ consists of a single point.  Since $u_\tau\in X_\tau$ always, we have $X_\tau=\{u_\tau\}$.

Consider any sequence $\la_n$ that tends nt to $\tau$.  Since $\tau$ is a $B$-point, $u_{\la_n}$ is bounded, by Corollary \ref{2.10}, and therefore has a weakly convergent subsequence, whose limit belongs to $X_\tau$, and hence is $u_\tau$.  As in the previous case it follows that $u_\la$ extends continuously to $S^-$ for any set $S$ that approaches $\tau$ nontangentially.  
\black
\end{proof}  
We can conclude the proof of Theorem \ref{4.19}.  Suppose (3): $\ph$ has a holomorphic differential at $\tau$ on $S$ and $|\ph(\tau)|=1$.  If $\tau\in\T^2$ then, by Lemma   \ref{constxdel}, $x_\tau$ is constant on $\Htau$,   and by Lemma \ref{xconstCpt} $u_\la\to u_\tau$ as $\la\to \tau$ in $S$.    On the other hand, if $\tau\in\T\times\D$, then by Lemma \ref{CptTD} the same conclusion follows.  In either case 
we deduce that (3) $\Rightarrow$ (1).
\end{proof}
Theorem \ref{4.19}  shows that, for a particular set $S$ untapered at $\tau$, the function-theoretic property of admitting a holomorphic differential on $S$ corresponds  to a continuity property of models of the function.  It leaves open the possibility that different untapered $S$ might behave differently, but in fact the 
\blue
following corollary essentially shows that the properties described are independent of the choice of $S$. \black
\begin{corollary}\label{4equiv}
Let  $\phi\in\Schur$, $\tau\in \tb$. 
 The following  six  conditions are equivalent. 
\begin{enumerate}
\item[\rm(1)] For every model  $(\mathcal{M}, u)$ of $ \phi$ and for every set $S \subset \mathbb{D}^2$ \blue that is untapered at $\tau$, $\lim_{\la \to\tau,\la\in S} |\phi(\la)|=1$ and \black  $u_\lambda$ extends  continuously  to $S^-$; 

\item[\rm(2)] there exist a model  $(\mathcal{M}, u)$ of $ \phi$ and a set $S \subset \mathbb{D}^2$ untapered at  $\tau$  \blue such that   $\lim_{\la \to\tau,\la\in S} |\phi(\la)|=1$ \black and  
$u_\lambda$ extends  continuously  to $S^-$; 

\item[\rm(3)] $\tau$ is a $C$-point for $\ph$, that is, for every set  $S \subset \mathbb{D}^2$   that approaches  $\tau$ nontangentially, $ \phi $ has  a  holomorphic  differential  at $ \tau$ on $ S$  with  $|\ph(\tau)| = 1$;

\item[\rm(4)] there exists a set  $S \subset \mathbb{D}^2$  untapered at  $\tau$ such that $ \phi $ has  a  holomorphic  differential  at $ \tau$ on $ S$  with  $|\ph(\tau)| = 1$; 

\item[\rm(5)]\blue $\tau$ is a B-point for $\ph$ and, \black for every model $(\M,u)$ of $\ph$, the function $x_\tau(.)$ defined by equation \eqref{defxtau} is constant on $\Htau$;

\item[\rm(6)]\blue $\tau$ is a B-point for $\ph$ and \black there exists a model $(\M,u)$ of $\ph$ such that $x_\tau(.)$ is constant on $\Htau$. 

\end{enumerate} 
\end{corollary}
\begin{proof} 
Trivially (1)$\Rightarrow$(2), 
(3)$\Rightarrow$(4) and (5)$\Rightarrow$(6).
By Theorem \ref{4.19},  (1)$\Leftrightarrow$(3) and (2)$\Leftrightarrow$(4).  Lemma \ref{constxdel} asserts that  (4)$\Rightarrow$(5). Hence also (2)$\Rightarrow$(5).

We claim that \blue (5)$\Rightarrow$(1) and (6)$\Rightarrow$(2).  Indeed, by Remark \ref{limit-mod}, the fact that $\tau$ is a B-point for $\ph$ implies that  $\lim_{\la \to\tau,\la\in S} |\phi(\la)|=1$. (5)$\Rightarrow$(1) then \black follows from Lemma \ref{xconstCpt} (for $\tau\in\T^2$) and Lemma \ref{CptTD} (for $\tau\in\T\times\D$).  Similarly (6)$\Rightarrow$(2).

We therefore have the implications
\[
\begin{array}{ccc} (3) & \Rightarrow & (4) \\
    \Updownarrow & &  \Updownarrow  \\
     (1) & \Rightarrow & (2)  \\
    \Uparrow & \raisebox{-0.7ex}{\rotatebox{45}{$\Leftarrow$}} & \Uparrow \\
     (5) & \Rightarrow & (6)
\end{array}
\]
from which the equivalence follows.
\end{proof} \black
\begin{remark}  \label{clarify}
$\tau\in\tb$ is a $C$-point for $\ph\in\Schur$ if and only if $\ph$ has a holomorphic differential on some set $S$ untapered at $\tau$ with $|\ph(\tau)|=1$.\\
 \rm  This fact was alluded to immediately after Definition \ref{4.4} and is simply the equivalence (3)$\Leftrightarrow$(4) of Corollary \ref{4equiv}.
\end{remark}
\begin{corollary}   \label{nabphtau}
If $\tau\in\tb$ is a $C$-point for $\ph\in\Schur$ and $(\M,u)$ is a model of $\ph$ then the angular gradient of $\ph$ at $\tau$ is given by
\beq \label{formanggrad}
\nabla \ph(\tau) = \ph(\tau)\begin{pmatrix}\overline{\tau^1}|| u^1_\tau ||^2 \\
   \overline{\tau^2}|| u^2_\tau ||^2 \end{pmatrix}.
\eeq
\end{corollary}
For Lemma \ref{constxdel} applies, and so the angular gradient is given by equation (\ref{diffonS}).

Note that, by Theorem \ref{2.16}, the sum of the moduli of the components of $\nabla\ph(\tau)$ is the lim inf in equation (\ref{normutau}).

\section{A Carath\'eodory theorem}
\label{carathm}
In this section we shall prove the second main result of the paper, Theorem \ref{5.3p}, an analog for the bidisk of Carath\'eodory's result that (C) $\Rightarrow$ (D) in Theorem \ref{thma1} and conversely.

In this (and only this) section we denote the directional derivative of  a holomorphic function $\phi$ in a direction $\delta \in \mathbb{C}^2$  by 
\newcommand\nad{\nabla_\delta}
$\nad\ph$.  We use this notation rather than the traditional $D_\delta \phi$  (as in Section \ref{directional}) to avoid confusion due to a surfeit of $D$s.

If $S$ is any untapered set at $\tau\in\tb$  
on which $\phi$ admits a holomorphic differential, then by Lemma \ref{constxdel}, 
\begin{equation} \label{5.2}
\phi(\tau+t \delta) = \omega + t(\delta^1 \omega \overline{\tau^1}||u_\tau^1||^2 + \delta^2 \omega \overline{\tau^2}||u_\tau^2||^2 ) + o(t),
 \end{equation}
and so
\beq\label{5.2a}
\nad\ph(\tau) = \omega \lan \de u_\tau, \tau u_\tau \ran.
\eeq
  This formula holds for all $\delta$ for which there exists $\epsilon >0$ such that $\tau + t\delta \in S$ for $t \in (0,\epsilon)$, that is, in the terminology of Section \ref{holdiffC}, for all $\de$ such that $S$ admits $-\de$. 

For Schur class functions of one variable, the \cjw Theorem tells us, {\em inter alia}, 
that if $\tau$ is a $B$-point of $f$, then (C) $f$ has a holomorphic differential on sets untapered at $\tau$, and (D) $f'(\la)$ converges to 
the angular derivative of $f$ at $\tau$ as $\la \stackrel{\rm nt}{\to} \tau$. As we observed in the Introduction,
the most straightforward analog of this statement for the Schur class in two variables fails. In the example in Section \ref{example}, $(1,1)$ is a $B$-point for $\ph$, but not a $C$-point, 
that is, $\ph$ does not have a holomorphic differential at $(1,1)$ on any untapered set.  
However, it {\em is } true in two variables that (C) implies (D): if $\tau$ is a $C$-point of $\ph$ then the gradient of $\ph$ at $\la$ tends to the angular gradient of $\ph$ at $\tau$ as $\la \stackrel{\rm nt}{\to} \tau$.  
Here, as usual, the gradient of $\ph$ at $\la$ means the vector
\[
\nabla \phi(\lambda) \df  \begin{pmatrix} \ds \frac{\partial \phi}{\partial \lambda^1}(\la) \\ \, \\  \ds\frac{\partial \phi}{\partial \lambda^2}(\la) \end{pmatrix}
\]
for any $\la \in \D^2$.
\black

\noindent {\em Proof of Theorem} \ref{5.3p}.  Suppose that $\tau\in\tb$ is a $C$-point for $\phi \in \mathcal{S}$.   Then $\ph(\tau)=\omega$ for some $\omega\in\T$. Let $S\subset \mathbb{D}^2$ be an untapered set at $\tau$, so that $S$ admits some open set $\Delta\subset \C^2$.   To prove that $\nabla\ph(\la) \to \nabla\ph(\tau)$ as $\la\to \tau$ in $S$ it will suffice to show that, for an open set of $\delta \in \mathbb{C}^2$,  
 $\nabla_\delta \phi(\lambda) \to \de\cdot\nabla  \ph(\tau)$  when $\la\to\tau, \,  \la\in S$.

Fix a model $(\mathcal{M},u)$ of $\phi$ and a realization $(a,\beta,\ga,D)$ of $(\M,u)$.    By Proposition \ref{getx}, for $\la\in\D^2$,
\[
1-\bar\omega \ph(\la) = \lan (1-\tau^*\la)u_\la, u_\tau \ran.
\]
Application of $ -\omega\nad$ to this equation yields: for all $\la\in\D^2,\, \de\in -\Htau$,
\begin{align*}
 \nad\ph(\la)& = \omega \lan \tau^*\de u_\la, u_\tau \ran - \omega\lan (1-\tau^*\la) \nad u_\la, u_\tau \ran\\   &= \omega \lan \de u_\la, \tau u_\tau \ran - \omega\lan (\tau-\la) \nad u_\la, \tau u_\tau \ran
\end{align*}
(for $\tau\in\T\times\D$ the last equation holds because $u_\tau^2=0$).
By equation (\ref{5.2a}), if $S$ admits $-\de$,
\[
\nad\ph(\tau) = \omega \lan \de u_\tau, \tau u_\tau \ran,
\]
and since $S$ admits $\Delta$, this equation holds in particular when $-\de\in \Delta$. 
By Corollary \ref{4equiv} $u_\la \to u_\tau$ for $\la\in S$ and we have $\lan \de u_\la, \tau u_\tau \ran \to \lan \de u_\tau, \tau u_\tau \ran$.   To  show that, for  $\de\in -\Delta$,  $\nad\ph(\la) \to \nad \ph(\tau)$ as $\la\to \tau, \, \la\in S$  it suffices to prove that, for all $\de \in -\Delta$, 
 \begin{equation} \label{5.4}
\lan(\tau-\lambda)   \nad u_\lambda, \tau u_\tau\ran \to 0 \quad \mathrm{as \ } \lambda\rightarrow \tau \mathrm{ \ in \ } S. 
\end{equation}
In fact we shall prove the stronger assertion, that 
\beq\label{stronger}
(\tau-\lambda)   \nad u_\lambda \to 0 \mbox{  as  } \la\to \tau, \quad \la\in S
\eeq
for every $\de \in -\Delta$.

   From the relation
\[
(1-D\lambda)u_\lambda = \gamma = (1-D\tau)u_\tau = (1-D\lambda)u_\tau +D(\lambda-\tau)u_\tau
\]
 we see that 
\begin{equation}\label{5.5}
 u_\lambda = u_\tau + (1-D\lambda)^{-1}D(\lambda-\tau) u_\tau.
\end{equation}
 
 This formula has two consequences of interest here. First, since Proposition \ref{4.19} guarantees that $u_\lambda$ extends to be continuous on $S^-$,
 \begin{equation} \label{5.6} 
(1-D\lambda)^{-1}D(\lambda-\tau)u_\tau \rightarrow 0 \quad \mathrm{as \ } \lambda\rightarrow \tau \mathrm{ \ in \ } S.
\end{equation}
 
 For the second consequence let us apply $\nad$ to equation (\ref{5.5}).  We have
\begin{align*}
\nad(1-D\la)^{-1} &= \lim_{h\to 0} \frac{1}{h}[(1-D(\la+h\de))^{-1} - (1-D\la)^{-1}]\\
   &=\lim_{h\to 0} \frac{1}{h}(1-D(\la+h\de))^{-1} [1-D\la - (1-D(\la+h\de))](1-D\la)^{-1}\\
  &= (1-D\la)^{-1}D\de (1-D\la)^{-1}.
\end{align*}
Hence, for all $\la\in\D^2$ and $\de\in -\Htau$,
 \begin{equation} \label{5.7}
 \nad u_\lambda = (1-D\lambda)^{-1}  D\de(1-D\lambda)^{-1}D{(\lambda-\tau) }u_\tau +(1-D\lambda)^{-1}  
 D\de u_\tau.
 \end{equation}

 We now turn to the proof of the limit relation (\ref{stronger}). In light of (\ref{5.7}) it suffices to establish the following two facts:  for every $\de\in -\Delta$,
\begin{equation} \label{5.8}
(\tau - \lambda)(1 - D\lambda)^{-1} D\de (1-D\lambda)^{-1}D(\lambda-\tau) u_\tau \rightarrow 0 \mathrm{\ as \ } \lambda\rightarrow \tau \mathrm{ \ in \ } S,
\end{equation} 
 \begin{equation} \label{5.9}
(\tau-\lambda)(1-D\lambda)^{-1}D\de u_\tau \rightarrow 0 \mathrm{\ as \ } \lambda\rightarrow \tau \mathrm{ \ in \ } S.  \end{equation}
 
 To see that (\ref{5.8}) holds note that since $S$ approaches $\tau$ nontangentially there exists a $c$ such that
 \[ 
||\tau - \lambda|| \le c (1 - ||\lambda||).
\]
 for all $\lambda \in S$.  Hence, for such $\lambda$,
 \begin{align*} ||(\tau -\lambda) (1-D\lambda)^{-1}D ||  & \le  ||\tau - \lambda||\cdot||(1-D\lambda)^{-1} D|| \\
 & \le c ( 1 -||\lambda||) \frac{1}{1 - ||\lambda||} \\
 & =c, 
\end{align*}
  that is,
 \begin{equation} \label{5.10} 
(\tau-\lambda)(1 - D\lambda)^{-1} D \quad\mathrm{\ is \ bounded \ on \ } S. 
\end{equation} 
The limit relation (\ref{5.8}) follows immediately by combination of this fact with (\ref{5.6}). 
 
 To see that (\ref{5.9}) holds, let $\de\in -\Delta$ and consider any sequence $\{\la_n\}\subset S$ converging to $\tau$.    Let $t_n= ||\la_n-\tau||$  and $\gamma_n =(\la_n-\tau)/t_n$, so that 
\[
\la_n = \tau +  t_n \ga_n, \quad ||\ga_n|| =1 \mbox{ and } t_n \to 0.
\]
 Define $\mu_n$ by $\mu_n = \tau + t_n \delta$. Since $\delta\in-\Delta$ and $t_n \to 0$ we have $\mu_n \in S$ for sufficiently large $n$, and further, $\mu_n \to \tau$.   By relation (\ref{5.6}),
\beq\label{limmun}
(1-D\mu_n)^{-1}D(\mu_n -\tau)u_\tau \to 0 \mbox{  as  } n \to \infty.
\eeq
We have 
 \begin{align}\label{5.10a} 
(\tau-\lambda_n)(1-&D{\lambda_n})^{-1}   D\de u_\tau  \nn\\
 & = (\tau - \lambda_n) [(1-D{\lambda_n})^{-1} - (1- D{\mu_n})^{-1}] D\de u_\tau\nn \\
 & \quad+ ( (\tau - \lambda_n) - (\tau - \mu_n)) (1-D{\mu_n})^{-1}   D\de u_\tau\nn \\
 &\quad + (\tau - \mu_n)(1 - D{\mu_n})^{-1}   D\de u_\tau 
\end{align}
 We verify (\ref{5.9}) by showing in succession that each of the three terms on the right hand side of equation (\ref{5.10a}) tends to $0$.  First observe that, since $\la_n-\mu_n = t_n(\ga_n-\de)$ and $\mu_n-\tau=t_n \de$,
\begin{align*} 
(\tau - \lambda_n)[(1-D{\lambda_n})^{-1}& -(1-D{\mu_n})^{-1}]D{\delta}u_{\tau}  \\
& = (\tau-\la_n) ( 1 -D{\lambda_n})^{-1}Dt_n(\ga_n-\delta)(1-D{\mu_n})^{-1}  D\de u_\tau \\
& =(\tau - \lambda_n)(1-D{\lambda_n})^{-1}D(\ga_n-\delta)(1-D{\mu_n})^{-1} D(\mu_n-\tau)u_\tau \\
& \rightarrow 0 
\end{align*}
since (\ref{5.10}) implies that $(\tau - \lambda_n)(1-D{\lambda_n})^{-1}D$ is bounded on $S$,  $||\ga_n-\de|| \leq 1+||\de||$ and  $(1 -D\mu_n)^{-1}D(\tau-\mu_n)u_\tau \rightarrow 0$ by relation (\ref{limmun}).

Next observe that 
\begin{align*} ((\tau - \lambda_n) - (\tau - \mu_n))(1-D{\mu_n})^{-1} D\de \mu_\tau & = t_n ( \delta-\ga_n)(1 - D{\mu_n})^{-1}  D\de u_\tau \\
& = (\gamma_n-\delta)(1-D{\mu_n})^{-1}D(\tau-\mu_n) u_\tau \\
& \rightarrow 0 
\end{align*}
by (\ref{limmun}). Finally, observe that
\begin{align*}(\tau - \mu_n)(1 - D{\mu_n})^{-1}  D\de u_\tau &  = -t_n \delta (1- D{\mu_n})^{-1} D\de u_\tau \\
& = \delta (1 - D{\mu_n})^{-1} D(\tau-\mu_n) u_\tau \\
& \rightarrow 0 
\end{align*}
by (\ref{limmun}). This establishes (\ref{5.9}) for every $\de\in-\Delta$ and so concludes the proof that $\nabla\ph(\la) \to \nabla\ph(\tau)$ as $\la\stackrel{\rm nt}{\to} \tau$.

Conversely, suppose that $\tau\in\tb$ is a $B$-point for $\ph$ such that $\nabla\ph(\la)\to\eta$ as $\la\stackrel{\rm nt}{\to} \tau$.  Consider a set $S$ untapered at $\tau$, so that $S$ admits some non-empty open set $\Delta$.  Since $\tau$ is a $B$-point $\ph$ is continuous on the closed line segment $[\tau,\la]$ joining $\tau$ to any point $\la\in\D^2$, and hence
\[
\ph(\la)=\ph(\tau) + \int_{[\tau,\la]} \nabla\ph(\zeta)\cdot {\rm d}\zeta.
\]
Hence, for any $\la\in\D^2$, the error
\begin{align*}
e(\la)  &\df \ph(\la)-\ph(\tau)- \eta\cdot (\la-\tau) \\
   &= \int_{[\tau,\la]} (\nabla\ph(\zeta)-\eta)\cdot {\rm d}\zeta.
\end{align*}
Let
\[
S' = S \cap \{\tau+t\de: \de\in\Delta, t>0\}.
\]
Then $S'$ is untapered at $\tau$, and $[\tau,\la] \subset S'$ for all $\la\in S'$ close enough to $\tau$.  Thus, for all such $\la$,
\[
\|e(\la)\| \leq (\sup_{\zeta\in[\tau,\la]} \|\nabla\ph(\zeta) - \eta\|) \: \|\la-\tau\|.
\]
Since $S'$ approaches $\tau$ nontangentially, the supremum on the right hand side tends to zero by hypothesis as $\la\to\tau$ in $S'$, and so $\|e(\la)\|=o(\|\la-\tau\|)$ as $\la\to\tau, \, \la\in S'$.  Thus $\ph$ has a holomorphic differential at $\tau$ on the untapered set $S'$, and so $\tau$ is a $C$-point for $\ph$, by Corollary \ref{4equiv}.
 \hfill $\Box$

\section{Concluding remarks} \label{conclud}
 Jafari \cite{jaf93} and Abate \cite{ab98} proved numerous results about limits of a Schur class function  and its difference quotients as the variable approaches a $B$-point along various special types of curve in a polydisk.  We summarise some of their results in our notation and terminology.     Despite the similarity of our titles, it will be clear that the focus of their work, on curves that approach $\tau$ in a tightly controlled way, is different from that of the present paper.

Consider $\ph\in\Schur$.  For simplicity we shall only consider $B$-points $\tau\in\tb$ (though Abate treats general $B$-points in $\partial(\D^n)$).
The classical Julia-Carath\'eodory Theorem can be applied to the function $z\mapsto \ph(z\tau)$ to obtain statements about limits along curves tending to $\tau$ and lying in the ``slice" $\D^2\cap \C\tau$.  The two authors show that similar statements hold for curves which, while not lying in this slice, approach it sufficiently closely as the curve tends to $\tau$.  Specifically, Theorems 0.3 and 3.1 of \cite{ab98} contain the following.

 {\em Suppose that
\[
\liminf_{\la\to\tau}\frac{1-|\ph(\la)|}{1-\|\la\|} =\al < \infty.
\]
Then there exists $\omega\in\T$ such that the Julia inequality {\rm (\ref{includ})} holds, so that $\ph(\la)\to\omega$ as $\la\to\tau$ horospherically.  Furthermore,
\begin{enumerate}
\item[\rm(i)] the difference quotient $(\ph(\la)-\omega)/(\tfrac 12 \lan \la,\tau\ran-1) $ has restricted $K$-limit $\al\omega$ at $\tau$;
\item[\rm (ii)] the difference quotient $(\ph(\la)-\omega)/(\la^j-\tau^j)$ has restricted $K$-limit $\al\omega\overline{\tau^j}$ at $\tau$;
\item[\rm(iii)] the directional derivative $D_\tau \ph(\la) $ has restricted $K$-limit $\al\omega$ at $\tau$;
\item[\rm(v)] $\partial\ph/\partial \la^j$ has a restricted $K$-limit at $\tau$, which is in general different from the corresponding incremental ratio.
\end{enumerate} }
In (i) and below the inner product $\lan \la,\tau\ran$ is with respect to the standard Euclidean structure on $\C^2$. 
We must explain restricted $K$-limits.   For $\tau\in\T^2$ we denote by $P_\tau$ the orthogonal projection operator onto the plane $\C\tau$ in $\C^2$, so that
\[
P_\tau \la = \frac{\lan \la,\tau\ran}{\lan \tau,\tau\ran}\tau = \tfrac 12 \lan \la,\tau\ran  \tau \qquad \mbox{ for } \la\in\C^2.
\]
\newcommand{\lt}{\lambda_t}
A function $\psi$ on $\D^2$ is said to have restricted $K$-limit $L$ at $\tau\in\T^2$ if $\psi(\lt)\to L$ as $t\to 1-$ whenever the following four conditions hold:
\begin{itemize}
\item $\lt, \ 0\leq t< 1$, is a continuous curve in $\D^2$ such that $\lim_{t\to1-} \lt = \tau$;
\item $\lt$ lies in a {\em Koranyi region} at $\tau$, that is
\[
\frac{|(\lt)^j - \tau^j|^2}{(1-|(\lt)^j|)(1-\|\lt\|)} \mbox{ is bounded for } 0 \leq t < 1, \ j=1,2;
\]
\item $\lt$ is {\em special}, that is, tangent to $\C\tau$ at $\tau$, meaning
\[
\frac{\lt - P_\tau\lt}{1-\|P_\tau\lt\|} \to 0 \mbox{  as  } t\to 1-;
\]
\item $\lt$ is {\em restricted}, that is, $P_\tau\lt \nt \tau$ as $t\to 1-$.
\end{itemize}
Lying in a Koranyi region at $\tau$ is intermediate between approaching nontangentially and plurinontangentially.  If $S=\{\lt : 0\leq t<1\}$ then
\[
S\mbox{ approaches }\tau \mbox{ nt } \Rightarrow S\mbox{ approaches $\tau$ in a Koranyi region }\Rightarrow S\mbox{ approaches $\tau$ pnt}.
\]

Jafari obtains results along similar lines to Abate. 
We should mention that the two authors also study holomorphic maps from a polydisk to another polydisk and to more general codomains.

Finally, we observe that we do not know whether our results extend to functions in the Schur class (rather than the Schur-Agler class) of the tridisk.

\bibliography{references}

\end{document}